\documentclass[fleqn,11pt,a4paper,leqno]{article}

\usepackage{graphicx}
\usepackage{graphics}
\usepackage{latexsym}
\usepackage{amssymb,amsmath}
\usepackage{multicol}
\usepackage{bussproofs}
\usepackage{enumerate}
\usepackage{xcolor}
\usepackage[latin1]{inputenc}
\usepackage{verbatim}

\title{}

\date{}

\numberwithin{equation}{section}

\newtheorem{definicion}{Definition}[section]

\newtheorem{definition}[definicion]{Definition}
\newtheorem{Lemma}[definicion]{Lemma}

\newtheorem{Theorem}[definicion]{Theorem}
\newtheorem{theorem}[definicion]{Theorem}

\newtheorem{lema}[definicion]{Lemma}

\newtheorem{teorema}[definicion]{Theorem}
\newtheorem{corolario}[definicion]{Corollary}

\newtheorem{remark}[definicion]{Remark}

\newenvironment{Proof}{\noindent\bf Proof \rm}{$\hfill
\square$}

\begin{document}

\title{Order in Implication Zroupoids}

\author{Juan M. CORNEJO and Hanamantagouda P. SANKAPPANAVAR}

\numberwithin{equation}{section}

\maketitle

\begin{abstract}
	
	 The variety $\mathbf{I}$ of implication zroupoids (using a binary operation $\to$ and a constant $0$) was defined and investigated by Sankappanavar in \cite{sankappanavarMorgan2012}, as a generalization of De Morgan algebras.  Also, in \cite{sankappanavarMorgan2012}, several new subvarieties of $\mathbf{I}$ were introduced, including the subvariety $\mathbf{I_{2,0}}$, defined by the identity: $x'' \approx x$, which plays a crucial role in this paper.  Some more new subvarieties of $\mathbf{I}$ are studied in  \cite{CoSa2015aI} that includes the subvariety $\mathbf{SL}$ of semilattices with a least element $0$; and an explicit description of semisimple subvarieties of $\mathbf{I}$ is given in \cite{CoSa2015semisimple}.

It is a well known fact that there is a partial order (denote it by  $\sqsubseteq$) induced by the operation $\land$, both in the variety $\mathbf{SL}$ of semilattices with a least element  and in the variety $\mathbf{DM}$ of De Morgan algebras.  As both $\mathbf{SL}$ and $\mathbf{DM}$ are subvarieties of $\mathbf{I}$ and the definition of partial order can be expressed in terms of the implication and the constant, it is but natural to ask whether the relation $\sqsubseteq$ on $\mathbf{I}$ is actually a partial order in some (larger) subvariety of $\mathbf{I}$ that includes both $\mathbf{SL}$ and $\mathbf{DM}$.

The purpose of the present paper is two-fold:  Firstly, a complete answer is given to the above mentioned problem.  Indeed, our first main theorem shows that the variety $\mathbf{I_{2,0}}$ is a maximal subvariety of $\mathbf{I}$ with respect to the property that the relation $\sqsubseteq$ is a partial order on its members..  In view of this result, one is then naturally led to consider the problem of determining the number of non-isomorphic algebras in $\mathbf{I_{2,0}}$ that can be defined on an $n$-element chain (herein called $\mathbf{I_{2,0}}$-chains), $n$ being a natural number.  Secondly, we answer this problem in our second main theorem which says that, for each $n \in \mathbb{N}$, there are exactly $n$ nonisomorphic $\mathbf{I_{2,0}}$-chains of size $n$.  

\end{abstract}

\thispagestyle{empty}

\section{Introduction}
%

\hspace{10pt} The widely known fact that Boolean algebras can be defined using only implication and a constant was extended to De Morgan algebras in \cite{sankappanavarMorgan2012}.  The crucial role played by a certain identity, called (I), led Sankappanavar, in 2012, to define and investigate, the variety $\mathbf{I}$ of implication zroupoids (I-zroupoids) generalizing De Morgan algebras.  Also, in \cite{sankappanavarMorgan2012}, he introduced several new subvarieties of $\mathbf{I}$ and found some relationships among those subvarieties.  One of the subvarieties of $\mathbf{I}$, called $\mathbf{I_{2,0}}$, defined by the identity: $x'' \approx x$ and studied in  \cite{sankappanavarMorgan2012}, plays a crucial role in this paper. 
In  \cite{CoSa2015aI}, we introduce several more new subvarieties of $\mathbf{I}$, including the subvariety $\mathbf{SL}$ which is term-equivalent to the (well known) variety of $\lor$-semilattices with a least element $0$, and describe further relationships among the subvarities of $\mathbf{I}$.  An explicit description of semisimple subvarieties of $\mathbf{I}$ is given in \cite{CoSa2015semisimple}. 

It is also a well known fact that there is a partial order induced by the operation $\land$, both in the variety $\mathbf{SL}$ of semilattices with a least element  and in the variety $\mathbf{DM}$ of De Morgan algebras.  As both $\mathbf{SL}$ and $\mathbf{DM}$ are subvarieties of $\mathbf{I}$ and the defintion of partial order can be expressed in terms of the implication and constant, it is but natural to ask whether the relation $\sqsubseteq$ (now defined) on $\mathbf{I}$ is actually a partial order in some (larger) subvariety of $\mathbf{I}$ that includes both $\mathbf{SL}$ and $\mathbf{DM}$.
   
The purpose of the present paper is two-fold:  Firstly, a complete answer is given to the above mentioned problem. 
Indeed, our first main theorem shows that the variety $\mathbf{I_{2,0}}$ 
 is a maximal subvariety of $\mathbf{I}$ with respect to the property that the relation $\sqsubseteq$, defined by: 
\begin{center}
	$x \sqsubseteq y$ if and only if $(x \to y')' =x$, for $x,y \in \mathbf{A}$ and $\mathbf{A} \in \mathbf{I}$,
\end{center}
is a partial order.  In view of this result, one is then naturally led to consider the problem of determining the number of non-isomorphic algebras in $\mathbf{I_{2,0}}$ ($\mathbf{I_{2,0}}$-chains) that can be defined on an $n$-element set, $n$ being a nutural number.  Secondly, we answer this problem in our second main result 
which says that, for each $n \in \mathbb{N}$, there are exactly $n$ nonisomorphic $\mathbf{I_{2,0}}$-chains of size $n$.   

\section{Preliminaries}

In this section we recall some definitions and results from \cite{CoSa2015aI}, {\cite{CoSa2015semisimple} and \cite{sankappanavarMorgan2012} that will be needed for this paper.  Basic references are \cite{balbesDistributive1974} and \cite{burrisCourse1981}.

\begin{definition} {\rm\cite{sankappanavarMorgan2012}}
	A {\em groupoid with zero} \rm{(}{\em zroupoid}, for short\rm{)} is an algebra $\mathbf A = \langle A, \to, 0\rangle$, where $\to$ is a binary operation and $0$ is a constant. A zroupoid $\mathbf A = \langle A, \to, 0\rangle$ is an {\em implication zroupoid} \rm{(}{\em \rm{I}-zroupoid}, for short) if the following identities hold in $\mathbf A$, where $x':= x \to 0$:
	\begin{itemize}
		\item[(I)] $(x \to y) \to z \approx [(z'\to x) \to (y \to z)']'$
		\item[${\rm(I_0)}$] $0''\approx 0$.
	\end{itemize}
\end{definition}

The variety of I-zroupoids is denoted by $\mathbf{I}$.

In this paper we use the characterizations of De Morgan algebras, Kleene algebras and Boolean algebras (see {\rm\cite{sankappanavarMorgan2012}}), and semilattices with least element $0$ (see {\rm \cite{CoSa2015aI}}), as definitions.

\begin{definition}
An implication zroupoid $\mathbf A = \langle A, \to, 0 \rangle$ is a {\it De Morgan algebra} {\rm(}$\mathbf{DM}$-algebra, for short{\rm)} if $\mathbf A$ satisfies the axiom:
\begin{itemize}
	\item
	[{\rm(DM)}]  $(x \to y) \to x \approx x. $
\end{itemize}

A $\mathbf{DM}$-algebra $\mathbf A = \langle A, \to, 0 \rangle$ is a {\it Kleene algebra} {\rm(}$\mathbf{KL}$-algebra, for short{\rm)} if
$\mathbf A$ satisfies the axiom:
\begin{itemize}
	\item[{\rm(KL$_1$)}] $(x \to x) \to (y \to y)' \approx x \to x$
\end{itemize}
or, equivalently,
\begin{itemize}
	\item[{\rm(KL$_2$)}] $(y \to y) \to (x \to x) \approx x \to x$.
\end{itemize}

A $\mathbf{DM}$-algebra $\mathbf A = \langle A, \to, 0 \rangle$ is a {\it Boolean algebra} {\rm(}$\mathbf{BA}$-algebra, for short{\rm)} if
$\mathbf A$ satisfies the axiom:
\begin{itemize}
	\item[{\rm(BA)}] $x \to x \approx 0'$.
\end{itemize}

An implication zroupoid $\mathbf A = \langle A, \to, 0 \rangle$ is a {\it semilattice with $0$} {\rm(}$\mathbf{SL}$-algebra, for short{\rm)} if $\mathbf A$ satisfies the axioms:
\begin{itemize}
	\item
	[{\rm(SM1)}]   $x'  \approx x$
         \item
	[{\rm(SM2)}]  $ x \to y \approx y \to x$.     {\rm (Commutativity)}.
\end{itemize}

We denote by $\mathbf{DM}$, $\mathbf{KL}$, $\mathbf{BA}$ and $\mathbf{SL}$, respectively, the variety of $\mathbf{DM}$-algebras,
$\mathbf{KL}$-algebras, $\mathbf{BA}$-algebras, and $\mathbf{SL}$-algebras. 

\end{definition}

We recall from \cite{sankappanavarMorgan2012} the definition of another subvariety of $\mathbf{I}$, namely 
$\mathbf I_{2,0}$, which plays a fundamental role in this paper.  
\begin{definition}
	 $\mathbf I_{2,0}$ denotes the subvariety of $\mathbf{I}$ defined by the identity:
	 $$x'' \approx x.$$
\end{definition}
We note that $\mathbf{DM}$, $\mathbf{KL}$, $\mathbf{BA}$ and $\mathbf{SL}$ are all subvarieties of $\mathbf I_{2,0}$ (see \cite{sankappanavarMorgan2012} and \cite{CoSa2015aI}).

\medskip

\begin{lema}{\rm\cite[Theorem 8.15]{sankappanavarMorgan2012}} \label{general_properties_equiv}
	Let $\mathbf A$ be an I-zroupoid. Then the following are equivalent:
	\begin{enumerate}[{\rm (a)}]
		\item $0' \to x \approx x$ \label{TXX} 
		\item $x'' \approx x$
		\item $(x \to x')' \approx x$ \label{reflexivity}
		\item $x' \to x \approx x$. \label{LeftImplicationwithtilde}
	\end{enumerate}
\end{lema}

\begin{lema}{\rm\cite{sankappanavarMorgan2012}} \label{general_properties}
	Let $\mathbf A \in \mathbf I_{2,0}$. Then
	\begin{enumerate}[{\rm (a)}]
		\item $x' \to 0' \approx 0 \to x$ \label{cuasiConmutativeOfImplic2}
		\item $0 \to x' \approx x \to 0'$. \label{cuasiConmutativeOfImplic}
	\end{enumerate}
\end{lema}

\medskip

Several identities true in ${\mathbf I_{2,0}}$ are given in \cite{CoSa2015aI}, \cite{CoSa2015semisimple} and \cite{sankappanavarMorgan2012}.  
Some of those that are needed for this paper are listed in the next lemma, which also presents some new identities of 
${\mathbf I_{2,0}}$ that will be useful later in this paper.  The proof of the lemma is given in the Appendix.

\begin{lema} \label{general_properties2}
	Let $\mathbf A \in \mathbf I_{2,0}$. Then $\mathbf A$ satisfies: 
	\begin{enumerate}[{\rm (1)}]
		\item $(x \to 0') \to y \approx (x \to y') \to y$ \label{281014_05} 
		\item $(0 \to x') \to (y \to x) \approx y \to x$ \label{281014_07} 
		\item $(y \to x)' \approx (0 \to x) \to (y \to x)'$ \label{291014_06}  
		\item $[x \to (y \to x)']' \approx (x \to y) \to x$ \label{291014_09}  
		\item $(y \to x) \to y \approx (0 \to x) \to y$ \label{291014_10} 
		\item $0 \to x \approx 0 \to (0 \to x)$ \label{311014_03}  
		\item $0 \to [(0 \to x) \to (0 \to y')'] \approx 0 \to (x \to y)$ \label{031114_05} 
		\item $[x' \to (0 \to y)]' \approx (0 \to x) \to (0 \to y)'$ \label{031114_06} 
		\item $0 \to (0 \to x)' \approx 0 \to x'$ \label{031114_07} 
		\item $0 \to (x' \to y)' \approx x \to (0 \to y')$ \label{071114_01} 
		\item $[(x \to 0') \to y]' \approx (0 \to x) \to y'$ \label{071114_02} 
		\item $0 \to [(0 \to x) \to y'] \approx x \to (0 \to y')$ \label{071114_03} 
		\item $0 \to (x \to y) \approx x \to (0 \to y)$ \label{071114_04} 
		\item $(x \to y) \to y' \approx y \to (x \to y)'$ \label{071114_05}  
		\item $(x' \to y) \to [(0 \to z) \to x'] \approx (0 \to y) \to [(0 \to z) \to x']$ \label{181114_07} 
		\item $0 \to (x \to y')' \approx 0 \to (x' \to y)$ \label{191114_05} 
		\item $x \to (y \to x') \approx y \to x'$ \label{281114_01} 
		%
		%
		\item $[(0 \to x) \to y] \to x \approx y \to x$ \label{291014_08}  
		\item $[0 \to (x \to y)] \to x \approx (0 \to y) \to x$ \label{311014_04} 
		\item $(0 \to x) \to (0 \to y) \approx x \to (0 \to y)$ \label{311014_06}  
		\item $x \to y \approx x \to (x \to y)$ \label{031114_04} 
		\item $[\{x \to (0 \to y)\} \to z]' \approx z \to [(x \to y) \to z]'$ \label{171114_01} 
		\item $[0 \to (x \to y)] \to y' \approx y \to (x \to y)'$ \label{181114_11} 
		\item $x \to [(y \to z) \to x]' \approx (0 \to y) \to [x \to (z \to x)']$ \label{181114_16} 
		\item $0 \to [(0 \to x) \to y] \approx x \to (0 \to y)$ \label{191114_02} 
		\item $x \to (y \to x)' \approx (y \to 0') \to x'$ \label{031214_16} 
		%
		%
		\item $[(x' \to y) \to (z \to x)'] \to [(y \to z) \to x] \approx (y \to z) \to x$ \label{281014_06}
		\item 
		     $[\{0 \to (x \to y)'\} \to (0 \to y')']' \approx 0 \to (x \to y)'$ \label{201114_02}
		\item 
		     $[[0 \to \{(x \to y) \to z\}] \to \{0 \to (y \to z)\}']' \approx 0 \to \{(x \to y) \to z\}$ \label{201114_03}
		\item $[x \to (0 \to y)']' \approx x' \to (y \to 0')'$ \label{201114_08}
		\item $[(0 \to x) \to y]' \approx y \to (x \to y)'$ \label{251114_03}
		\item $[x \to (y \to 0')']' \approx x' \to (0 \to y)'$ \label{271114_01}
		\item $(x \to y)' \to (0 \to x)' \approx y' \to x'$ \label{271114_02}
		\item $(0 \to x)' \to (0 \to y)'\approx 0 \to (x' \to y')$ \label{011214_01}
		\item  
		     $[(x \to y)' \to \{y \to (x \to y)'\}']' \approx (x \to y)'$ \label{011214_03}
		\item 
		    $[\{(0 \to x) \to y\} \to (x \to y)']' \approx (0 \to x) \to y$ \label{011214_04}
		\item 
		   $[\{x \to (y \to x)'\} \to x]' \approx x \to (y \to x)' $. \label{031214_07}
	\end{enumerate}
\end{lema}

\section{Partial order in Implication Zroupoids}

Let $\mathbf A = \langle A; \to, 0\rangle \in \mathbf I$.  We define the operations $\land$ and $\lor$ on $\mathbf A$ by: 

\begin{itemize}
	\item 	$x \land y := (x \to y')'$,  
	\item  $x \lor y := (x' \land y')'$.
\end{itemize}

Note that the above definition of $\land$ is a simultaneous generalization of the $\land$ operation of De Morgan algebras and that of $\mathbf{SL}$ (= semilattices with least element $0$).   It is, of course, well known that the meet operation induces a partial order on both $\mathbf{DM}$ and $\mathbf{SL}$,  which naturally leads us to the following definition of a binary relation $\sqsubseteq$ on algebras in $\mathbf{I}$.

\begin{definicion} \label{order_relation}
	Let $\mathbf A \in \mathbf I$. We define the relation $\sqsubseteq$ on $A$ as follows:
	$$x \sqsubseteq y \mbox{ if and only if } x \land y = x  \quad ( equivalently, (x \to y')' = x).$$
	For $a,b \in A$, we write 
	\begin{itemize}
		\item $a \sqsubset b$ if $a \sqsubseteq b$ and $a \not= b$,  
		\item $a \sqsupseteq b$ if $b \sqsubseteq a$, and
		\item $a \sqsupset b$ if $a \sqsupseteq b$ and $a \not= b$. 
	\end{itemize}
\end{definicion}

We already know from \cite{CoSa2015aI} that $\langle A; \land, \lor \rangle$ is a lattice if and only if $\mathbf A$ is a De Morgan Algebra, implying that $\sqsubseteq$ is a partial order on $A$.  We know (see \cite{CoSa2015aI}) that $\sqsubseteq$ is also a partial order on algebras in $\mathbf{SL}$.  This fact led us naturally to consider the possibility of the existence of a  subvariety $\mathbf{V}$ of $\mathbf{I}$, containing both $\mathbf{SL}$ and $\mathbf{DM}$, such that, for every algebra $\mathbf{A}$ in $\mathbf{V}$, the relation $\sqsubseteq$ on $\mathbf{A}$ is actually a partial order.


In this section we will prove our first main result which says that the subvariety  $\mathbf{I_{2,0}}$,  is a maximal
subvariety of $\mathbf{I}$ with respect to the property that the relation $\sqsubseteq$ 
is a partial order on every member of that variety.  To achieve this end, we need to, first, prove that $\sqsubseteq$ is indeed a partial order on every member of $\mathbf{I_{2,0}}$, which will be done using the following lemmas.

\begin{lema} \label{antisymmetry}
Let $\mathbf A \in \mathbf I_{2,0}$.  Then the relation $\sqsubseteq$ is antisymmetric on  $\mathbf A$.   
\end{lema}

\begin{Proof}
Let $a,b \in A$ such that $a \sqsubseteq b $ and $b \sqsubseteq a $.  
	Let $c \in A$ be arbitrary. Then, using (I) and the hypothesis, one observes that $(c \to a) \to b' = [(b \to c) \to (a \to b')']' = [(b \to c) \to a]'$. Consequently,
	\begin{equation}\label{281014_03}
	(c \to a) \to b' = [(b \to c) \to a]',  \mbox{ where } c \in A.
	\end{equation}
	Hence,
	$$
	\begin{array}{lcll}
	a'  & = & (a \land b)' & \mbox{by hypothesis}\\
	   & = & (a \to b')'' & \mbox{by definition of $\land$} \\
	  & = & a \to b' & \mbox{} \\
	 & = & (a' \to a) \to b'  & \mbox{using Lemma } \ref{general_properties_equiv} (\ref{LeftImplicationwithtilde})  \\
	& = & [(b \to a') \to a]' & \mbox{from (\ref{281014_03}) with } c=a' \\
	& = & [(b \to a')'' \to a]' & \mbox{} \\
	& = & (b' \to a)' & \mbox{by hypothesis},
	\end{array}
	$$
	and, therefore,
	\begin{equation}\label{281014_04}
	a' = (b' \to a)'.
	\end{equation}	
Now,  
	$$
	\begin{array}{lcll}
	b' & =  & [b \to a']'' & \mbox{by hypothesis} \\
	& = & b \to a' & \mbox{} \\
        & = & (0 \to a'') \to (b \to a') & \mbox{by Lemma \ref{general_properties2} (\ref{281014_07}) with }  x=a', y = b \\
	& = &(0 \to a) \to (b \to a')'' & \mbox{} \\
	& = & (0 \to a) \to b' & \mbox{by hypothesis.}
	\end{array}
	$$
	Thus,
	\begin{equation}\label{291014_01}
	b' = (0 \to a) \to b'.
	\end{equation}
	Therefore,
	$$
	\begin{array}{lcll}
	a'  & = & [b' \to a]' & \mbox{from (\ref{281014_04})} \\
	& = & [(b \to 0) \to a]' & \mbox{} \\
	& = & (0 \to a) \to b' &  \mbox{from (\ref{281014_03}) with } c=0 \\
	& = & b' & \mbox{by (\ref{291014_01})}.
	\end{array}
	$$
	Consequently, we have that $a = a'' = b'' = b$, thus proving that $\sqsubseteq$ is antisymmetric on $\mathbf A$.	
\end{Proof}
\medskip

Now, we turn to proving the transitivity of the relation $\sqsubseteq$.  For this, we need the following lemmas.  
The proof of the following (technical) lemma is given in the Appendix.

\begin{lema} \label{pretransitivity1}
Let $\mathbf A \in \mathbf I_{2,0}$ with $a,b \in A$ such that $a \sqsubseteq b$.  Let $d \in A$ be arbitrary.  Then 
	\begin{enumerate}[{\rm(1)}]
		\item $(0 \to a') \to b = a' \to b$ \label{031114_10}
		\item $b \to a' = (0 \to b) \to a'$ \label{031114_11}
		\item $b \to a' = a'$ \label{031114_12}
		\item $0 \to (a' \to b) = 0 \to a$ \label{201114_04}
		\item $[(b \to d) \to a]' = (d \to a) \to b'$ \label{211114_03}
		\item $(0 \to d) \to a' = [\{d \to (0 \to b')\} \to a]'$ \label{211114_01}
		\item $a \to [(a' \to d) \to \{(0 \to a) \to b'\}] = (0 \to d) \to a'$ \label{211114_02}
		\item $a \to [(d \to a) \to b'] = a \to (d \to a)'$ \label{211114_04}
		\item $[0 \to (b \to d)] \to a = (0 \to d) \to a$ \label{211114_06}
		\item $[b \to (a \to d)] \to a = (0 \to d) \to a$ \label{211114_07}
		\item $b \to (0 \to a') = 0 \to a'$ \label{031214_11}
		\item $[(d \to a) \to b']' = (b \to d) \to a$ \label{031214_13}
		\item $a' \to b = b' \to a$ \label{031214_14}
		\item $(d \to a') \to b = (d \to 0') \to (a' \to b)$ \label{031214_15}
		\item $[(0 \to a') \to b]' = (0 \to a) \to b'$ \label{031214_21}
		\item $(a' \to b)' = (0 \to a) \to b'$ \label{031214_17}
		\item $b' \to [(b \to d) \to a] \sqsubseteq 0 \to b$. \label{031214_19}
		
	\end{enumerate}
\end{lema}

\medskip


\begin{lema} \label{pretransitivity2}
	Let $\mathbf A \in \mathbf I_{2,0}$ and let $a,b, e \in A$ such that $(a \to b')' = a$ and $(0 \to e') \to b = b$, and let $d \in A$ be arbitrary.  Then
	\begin{enumerate}[(a)]
		\item $b \to d = (0 \to (d \to e)) \to (b \to d)$  \label{031114_14}
		\item $(0 \to e) \to a' = a'$ \label{181214_02}
		\item $(0 \to e') \to a = a$. \label{011214_10}
		\item $(0 \to e) \to [a \to (a \to d)] = a \to d$. \label{011214_09}
	\end{enumerate}
\end{lema}

\begin{Proof}
	\begin{itemize}
		
		\item[(\ref{031114_14})]
		$$
		\begin{array}{lcll}
		b \to d & = & [(0 \to e') \to b] \to d & \mbox{by hypothesis} \\
		& = & [\{d' \to (0 \to e')\} \to (b \to d)'] \to [\{(0 \to e') \to b\} \to d] & \mbox{by Lemma \ref{general_properties2} (\ref{281014_06}) } \\
		&   & \hspace{5cm} \mbox{using } x = d, y = 0 \to e', z = b \\
		& = & [\{d' \to (0 \to e')\} \to (b \to d)'] \to (b \to d) & \mbox{by hypothesis} \\
		& = & [\{d' \to (0 \to e')\} \to 0'] \to (b \to d) & \mbox{by Lemma \ref{general_properties2} (\ref{281014_05})} \\
		& = & [0 \to \{d' \to (0 \to e')\}'] \to (b \to d) & \mbox{by Lemma \ref{general_properties} (\ref{cuasiConmutativeOfImplic2})} \\
		& = & [0 \to \{(0 \to d) \to (0 \to e')'\}] \to (b \to d) & \mbox{by Lemma \ref{general_properties2} (\ref{031114_06})} \\
		& = & [0 \to (d \to e)] \to (b \to d) & \mbox{by Lemma \ref{general_properties2} (\ref{031114_05})}.
		\end{array}
		$$

		\item[(\ref{181214_02})]
		
		Using Lemma \ref{pretransitivity1} (\ref{031114_12}) (twice), and (\ref{031114_14}) with $d = a'$, we obtain $[0 \to (a' \to e)] \to a' = [0 \to (a' \to e)] \to (b \to a') = b \to a' = a'$. Hence,
		\begin{equation} \label{031114_15}
		[0 \to (a' \to e)] \to a' = a'.
		\end{equation}
		Then,
		$$
		\begin{array}{lcll}
		(0 \to e) \to a'  & = & [0 \to (a' \to e)] \to a' & \mbox{by Lemma \ref{general_properties2} (\ref{311014_04}) using } x = a', y = e \\
		& = & a' & \mbox{by (\ref{031114_15})}.
		\end{array}
		$$
		\item[(\ref{011214_10})]
		$$
		\begin{array}{lcll}
		(0 \to e') \to a & = & [0 \to (0 \to e)'] \to a & \mbox{by Lemma \ref{general_properties2} (\ref{031114_07})} \\
		& = & [(0 \to e) \to 0'] \to a & \mbox{by Lemma \ref{general_properties} (\ref{cuasiConmutativeOfImplic2})} \\
		& = & [(0 \to e) \to a'] \to a & \mbox{by Lemma \ref{general_properties2} (\ref{281014_05})} \\
		& = & a' \to a & \mbox{by (\ref{181214_02})}\\  
		& = & a & \mbox{by Lemma \ref{general_properties_equiv} (\ref{LeftImplicationwithtilde})}.
		\end{array}
		$$
		
		\item[(\ref{011214_09})]

		$$
		\begin{array}{lcll}
		a \to d & =  & [(0 \to e') \to a] \to d & \mbox{by item (\ref{011214_10})} \\
		& = & [\{d' \to (0 \to e')\} \to (a \to d)'] \to [\{(0 \to e') \to a\} \to d] & \mbox{by Lemma \ref{general_properties2} (\ref{281014_06}) with }  \\
		&  &  &  x = d, y = 0 \to e', z = a \\
		& = & [\{d' \to (0 \to e')\} \to (a \to d)'] \to (a \to d) &  \mbox{by item (\ref{011214_10})}  \\
		& = & [\{d' \to (0 \to e')\} \to 0'] \to (a \to d) & \mbox{by Lemma \ref{general_properties2} (\ref{281014_05})} \\
		& = & [0 \to \{d' \to (0 \to e')\}'] \to (a \to d)  & \mbox{by Lemma \ref{general_properties} (\ref{cuasiConmutativeOfImplic2})} \\
		& = & [0 \to \{(0 \to d) \to (0 \to e')'\}] \to (a \to d) & \mbox{by Lemma \ref{general_properties2} (\ref{031114_06}) with } \\
		&  &  & x = d, y = e' \\
		& = & [0 \to (d \to e)]  \to (a \to d) & \mbox{by Lemma \ref{general_properties2} (\ref{031114_05})}
		\end{array}
		$$
		Thus,
		\begin{equation} \label{031114_18}
		a \to d = [0 \to (d \to e)]  \to (a \to d).
		\end{equation}
		Now,


		\begin{align*} 
	          (0 \to e) \to [a \to (a \to d)]  
		& =  [0 \to [\{a \to (a \to d)\} \to e]] \to [a \to (a \to d)]\\
		 & \quad \quad \text{by Lemma \ref{general_properties2}  (\ref{311014_04}) with }
		            x = a \to (a \to d), y = e \\
		& =  [0 \to [\{a \to (a \to d)\} \to e]] \to [a \to \{a \to (a \to d)\}] \\
		 & \quad \quad \text{by Lemma \ref{general_properties2} (\ref{031114_04})} \\
		& =  a \to [a \to (a \to d)]\\
		 & \quad \quad \text{by (\ref{031114_18}) replacing } d \mbox{ with } a \to (a \to d) \\
                  & =  a \to d\\
                   & \quad \quad \text{by Lemma \ref{general_properties2} (\ref{031114_04})}.
		\end{align*}
	\end{itemize} 
	Thus, (d) is proved and the proof of the lemma is complete.	
\end{Proof}

\medskip

Each of the next three lemmas prove a crucial step in the proof of transitivity of $\sqsubseteq$. 

\begin{lema} \label{transitivity}  
	Let $\mathbf A \in \mathbf I_{2,0}$ and let $a,b,c \in A$ such that $a \sqsubseteq b $ and $b \sqsubseteq c$. Let $d,e,f  \in A$ be arbitrary.  Then 
\begin{enumerate}[{\rm(1)}]
\item
    $(0 \to c') \to b = b$ \label{031114_13} 
 
\item   
    $(0 \to c) \to [a \to (a \to d)] = a \to d$ \label{041114_01} 
 
\item  $(0 \to c) \to (a \to d) = a \to d$ \label{041114_02}

\item  $\left[0 \to \left((0 \to b) \to c' \right) \right] \to b = b$    \label{101114_01}

\item 
        $\left\{d' \to \left[0 \to ((0 \to b) \to c') \right] \right\} \to (b \to d)' = (b \to d)' $ \label{101114_02}
\item
       $(b \to d) \to [e \to (b \to d)]' = [e \to 0'] \to (b \to d)' $ \label{101114_03}

\item 
	$[b \to (a \to c')] \to a = a$ \label{101114_04}

\item  
	$(0 \to b) \to (a \to d) = a \to d$ \label{191114_01}

\item 
	$0 \to [b \to (a \to d)] = 0 \to (a \to d)$ \label{191114_03}

\item  
	$0 \to [\{b \to (a \to d)\} \to e] = 0 \to [(a \to d) \to e]$ \label{191114_04}

\item  
	$[0 \to (d' \to c)] \to (0 \to b)' = (0 \to d) \to (0 \to b)'$ \label{201114_01}

\item 
	$0 \to (a' \to c) \sqsubseteq 0 \to b$ \label{201114_05}

\item  
	$(0 \to a) \to (0 \to b)' = (0 \to a)' $ \label{201114_06}

\item 
	$0 \to (a' \to c) = 0 \to a$ \label{201114_07}

\item 
	$(d \to e) \to [\{b \to (a \to f)\}' \to (0 \to a)'] = (d \to e) \to [(a' \to b) \to (f' \to a')]$.  \label{201114_09}
\end{enumerate}
\end{lema}

\begin{Proof}  
By hypothesis, we have $(a \to b')' = a$ and $(b \to c')' = b$.\\	
         
 \begin{enumerate} [{\rm(1)}]
 \item    
	$$
	\begin{array}{lcll}
	(0 \to c') \to b & = & (c \to 0') \to b & \mbox{by Lemma \ref{general_properties} (\ref{cuasiConmutativeOfImplic2})} \\
	& = & [(b' \to c) \to (0' \to b)']' & \mbox{by (I)} \\
	& = & [(b' \to c) \to b']' & \mbox{by Lemmaa \ref{general_properties_equiv} (\ref{TXX})}\\
	& = & [(0 \to c) \to b']' & \mbox{by Lemma \ref{general_properties2} (\ref{291014_10})} \\
	& = & [(c' \to 0') \to b']' & \mbox{by Lemma \ref{general_properties} (\ref{cuasiConmutativeOfImplic2})} \\
	& = & [(b'' \to c') \to (0' \to b')']'' & \mbox{from (I)} \\
	& = & (b'' \to c') \to (0' \to b')' & \mbox{} \\
	& = & (b \to c') \to (0' \to b')' & \mbox{} \\
	& = & (b \to c') \to b'' & \mbox{by Lemmaa \ref{general_properties_equiv} (\ref{TXX})}\\
	& = & (b \to c') \to b & \mbox{} \\
	& = & b' \to b & \mbox{by hypothesis} \\
	& = & b & \mbox{by Lemma \ref{general_properties_equiv} (\ref{LeftImplicationwithtilde})}.
	\end{array}
	$$
	
 \item  
	This is immmediate from (1) and Lemma \ref{pretransitivity2} (\ref{011214_09}) with $e = c$.\\ 
	
	
\item  
	Using Lemma \ref{general_properties2} (\ref{031114_04}) and %
	 (\ref{041114_01}) we have that $(0 \to c) \to (a \to d) = (0 \to c) \to [a \to (a \to d)] = a \to d$, implying (3). \\
\item  
	$$
	\begin{array}{lcll}
	\left[0 \to \left((0 \to b) \to c' \right) \right] \to b  & = & \left\{(b' \to 0) \to \left[\left((0 \to b) \to c' \right) \to b \right]' \right\}' & \mbox{by (I)} \\
	& = & \left\{b \to \left[\left((0 \to b) \to c' \right) \to b \right]' \right\}' & \mbox{} \\
	& = & \left\{b \to (c' \to b)' \right\}' & \mbox{by Lemma \ref{general_properties2} (\ref{291014_08})} \\
	& = & \left\{(b' \to b) \to (c' \to b)' \right\}' & \mbox{by Lemma \ref{general_properties_equiv} (\ref{LeftImplicationwithtilde})} \\
	& = & (b \to c') \to b & \mbox{by (I)} \\
	& = & (b \to c')'' \to b & \mbox{} \\
	& = & b' \to b & \mbox{by hypothesis} \\
	& = & b & \mbox{by Lemma \ref{general_properties_equiv} (\ref{LeftImplicationwithtilde})}.
	\end{array}
	$$
\item	 
	
	$$
	\begin{array}{lcll}
	\left\{d' \to \left[0 \to ((0 \to b) \to c') \right] \right\} \to (b \to d)' & = & \left\{\left[\left[0 \to \left((0 \to b) \to c' \right) \right] \to b \right] \to d \right\}' & \mbox{by (I)} \\
	& = & (b \to d)' & \mbox{by (\ref{101114_01})}.
	\end{array}
	$$
\item	 
	$$
	\begin{array}{lcll}
	(b \to d) \to [e \to (b \to d)]'  & = & [e \to (b \to d)] \to (b \to d)' & \mbox{by Lemma \ref{general_properties2} (\ref{071114_05}) with } \\
	&  &  & x = e, y = b \to d \\
	& = & [e \to 0'] \to (b \to d)' & \mbox{by Lemma \ref{general_properties2} (\ref{281014_05})}.
	\end{array}
	$$
\item 
	$$
	\begin{array}{lcll}
	[b \to (a \to c')] \to a & = & [(a' \to b) \to \{(a \to c') \to a\}']' & \mbox{by (I)} \\
	& = & [(a' \to b) \to \{(0 \to c') \to a\}']' & \mbox{by Lemma \ref{general_properties2} (\ref{291014_10})} \\
	& = & [(a' \to b) \to a']' & \mbox{by (\ref{031114_13}) and Lemma \ref{pretransitivity2} (\ref{011214_10})} \\
	& = & [(0 \to b) \to a']' & \mbox{by Lemma \ref{general_properties2} (\ref{291014_10})} \\
	& = & (a \to 0) \to (b \to a')' & \mbox{by (I)} \\
	& = & a' \to (b \to a')' & \mbox{} \\
	& = & a' \to a'' & \mbox{by Lemma \ref{pretransitivity1} (\ref{031114_12})} \\
	& = & a' \to a & \mbox{} \\
	& = & a & \mbox{by Lemma \ref{general_properties_equiv} (\ref{LeftImplicationwithtilde})}.
	\end{array}
	$$
\item 

	$$
	\begin{array}{lcll}
	(0 \to b') \to b & = & (0 \to 0') \to b & \mbox{by Lemma \ref{general_properties2} (\ref{281014_05})} \\
	& = & (0'' \to 0') \to b & \mbox{} \\
	& = & 0' \to b & \mbox{by Lemma \ref{general_properties_equiv} (\ref{LeftImplicationwithtilde})} \\
	& = & b & \mbox{by Lemma \ref{general_properties_equiv} (\ref{TXX})}.
	\end{array}
	$$
	Hence, by the hypothesis, together with Lemma \ref{pretransitivity2} (\ref{011214_09}),  we obtain that
	$(0 \to b) \to \{a \to (a \to d)\} = a \to d$. Hence, by Lemma \ref{general_properties2} (\ref{031114_04}), we have
	$(0 \to b) \to (a \to d) = a \to d.$
\item  
	$$
	\begin{array}{lcll}
	0 \to (a \to d) & = & 0 \to [(0 \to b) \to (a \to d)] & \mbox{by (\ref{191114_01})} \\
	& = & b \to [0 \to (a \to d)] & \mbox{by Lemma \ref{general_properties2} (\ref{191114_02}) with } x = b, y = a \to d \\
	& = & 0 \to [b \to (a \to d)]. & \mbox{by Lemma \ref{general_properties2} (\ref{071114_04})}.
	\end{array}
	$$
	
\item  
	$$
	\begin{array}{lcll}
	0 \to [\{b \to (a \to d)\} \to e] & = & [b \to (a \to d)] \to (0 \to e) & \mbox{by Lemma \ref{general_properties2} (\ref{071114_04})} \\
	& = & 0 \to [[0 \to \{b \to (a \to d)\}] \to e] & \mbox{by Lemma \ref{general_properties2} (\ref{191114_02})} \\
	& = & 0 \to [\{0 \to (a \to d)\} \to e] & \mbox{by (\ref{191114_03})} \\
	& = & (a \to d) \to (0 \to e) & \mbox{by Lemma \ref{general_properties2} (\ref{191114_02})} \\
	& = & 0 \to [(a \to d) \to e] & \mbox{by Lemma \ref{general_properties2} (\ref{071114_04})}.
	\end{array}
	$$
\item 
	$$
	\begin{array}{lcll}
	[0 \to (d' \to c)] \to (0 \to b)' & = & [0 \to (d' \to c)] \to (b' \to 0')' & \mbox{by Lemma \ref{general_properties} (\ref{cuasiConmutativeOfImplic2})} \\
	& = & [\{(d' \to c) \to b'\} \to 0']' & \mbox{by (I)} \\
	& = & [\{(b \to d') \to (c \to b')'\}' \to 0']' & \mbox{by (I)} \\
	& = & [\{(b \to d') \to b''\}' \to 0']' & \mbox{by Lemma \ref{pretransitivity1} (\ref{031114_12})} \\
	& = & [\{(b \to d') \to b\}' \to 0']' & \mbox{} \\
	& = & [\{(0 \to d') \to b\}' \to 0']' & \mbox{by Lemma \ref{general_properties2} (\ref{291014_10})} \\
	& = & [0 \to \{(0 \to d') \to b\}]' & \mbox{by Lemma \ref{general_properties} (\ref{cuasiConmutativeOfImplic2})} \\
	& = & [(0 \to d') \to (0 \to b)]' & \mbox{by Lemma \ref{general_properties2} (\ref{071114_04})} \\
	& = & [(d \to 0') \to (0 \to b)]' & \mbox{} \\ 
	& = & (0 \to d) \to (0 \to b)' & \mbox{by Lemma \ref{general_properties2} (\ref{071114_02})}.
	\end{array}
	$$
\item 
	$$
	\begin{array}{lcll}
	0 \to (a' \to c) & = & 0 \to [(a \to b')'' \to c] & \mbox{by hyphotesis} \\
	& = & 0 \to [(a \to b') \to c] & \mbox{} \\
	& \sqsubseteq  & 0 \to (b' \to c) & \mbox{by Lemma \ref{general_properties2} (\ref{201114_03})} \\
	& = & 0 \to b & \mbox{by hyphotesis and Lemma \ref{pretransitivity1} (\ref{201114_04})}.
	\end{array}
	$$
         
\item 
	$$
	\begin{array}{lcll}
	(0 \to a) \to (0 \to b)'  & = & [a' \to (0 \to b)]' & \mbox{by Lemma \ref{general_properties2} (\ref{031114_06})} \\
	& = & [0 \to (a' \to b)]' & \mbox{by Lemma \ref{general_properties2} (\ref{071114_04})} \\
	& = & (0 \to a)'. & \mbox{by hyphotesis and Lemma \ref{pretransitivity1} (\ref{201114_04})}.
	\end{array}
	$$
	
\item  
	$$
	\begin{array}{lcll}
	0 \to (a' \to c) & = & [\{0 \to (a' \to c)\} \to (0 \to b)']' & \mbox{by (\ref{201114_05})} \\
	& = & [(0 \to a) \to (0 \to b)']' & \mbox{by (\ref{201114_01}) with } d = a \\
	& = & (0 \to a)'' & \mbox{by (\ref{201114_06})} \\
	& = & 0 \to a. & \mbox{}
	\end{array}
	$$
\item 
	\begin{align*}
	(d \to e) \to [(a' \to b) \to (f' \to a')] \\
	 &\hspace{-2cm}=  (d \to e) \to [\{(0 \to a') \to b\} \to (f' \to a')] \\
	           & \quad \mbox{by Lemma \ref{pretransitivity1} (\ref{031114_10})} \\
	&\hspace{-2cm} =  (d \to e) \to [\{(0 \to a') \to b\} \to \{(f \to 0) \to a'\}] 
	                     & \mbox{} \\
	&\hspace{-2cm} =  (d \to e) \to [\{(0 \to a') \to b\} \to \{(a \to f) \to (0 \to a')'\}'] \\
	              &\quad \mbox{by (I)} \\
	&\hspace{-2cm} =  (d \to e) \to [\{b \to (a \to f)\} \to (0 \to a')']' \\
	             &\quad \mbox{by (I)} \\
	&\hspace{-2cm} =  (d \to e) \to [\{b \to (a \to f)\}' \to (a' \to 0')']  \\
	  & \quad \mbox{by (\ref{201114_08}) with } x = b \to (a \to f) \mbox{ and } y = a' & \\
	&\hspace{-2cm} =  (d \to e) \to [\{b \to (a \to f)\}' \to (0 \to a)']\\
	 &\quad \mbox{by Lemma \ref{general_properties} (\ref{cuasiConmutativeOfImplic2})}.
	\end{align*}
	Hence, we have
	$(d \to e) \to [\{b \to (a \to f)\}' \to (0 \to a)'] = (d \to e) \to [(a' \to b) \to (f' \to a')].$
\end{enumerate}	
\end{Proof}

\begin{Lemma} \label{Lemma_36}
Let $\mathbf A \in \mathbf I_{2,0}$ and let $a,b,c \in A$ such that $a \sqsubseteq b $ and $b \sqsubseteq c $. Let $d  \in A$ be arbitrary.  Then 
\begin{enumerate} [{\rm(a)}]
\item 
	$[c \to (b \to a')] \to b = (0 \to a') \to b$ \label{211114_08} 

\item 
	$(c \to a') \to b = a' \to b$ \label{211114_10}

\item 
	$(a' \to b) \to (c \to a') = c \to a'$ \label{251114_01}

\item 
	$c \to a' = a \to [b \to (a \to c')]$ \label{251114_02}

\item 
	$0 \to (a \to d) = 0 \to [c \to (a \to d)]$ \label{011214_02}

\item 
	$(d \to a) \to d \sqsubseteq (a' \to b) \to d$ \label{011214_05}

\item 
	$(a' \to b) \to c' = (0 \to a) \to b' $\label{011214_07} 

\item 
	$0 \to (a \to c') \sqsubseteq 0 \to a'  $ \label{011214_08} 

\item
	$0 \to (a \to c') = 0 \to a'.$ \label{021214_01} 

\item 
	$c \to (a \to c') \sqsubseteq 0 \to (a \to c')$ \label{021214_02}  

\item 
        $c \to (a \to c') \sqsubseteq 0 \to a'$ \label{031214_03}  


\item 
	$(c \to (a \to c'))' \to (0 \to a)' = c \to (a \to c')$  \label{031214_04} 

\item 
	$a \to [b \to (a \to c')] = a \to c' $ \label{031214_05}    
	
\item
	$c \to a' = a \to c'$. \label{031214_06}  

\end{enumerate}

\end{Lemma}

\begin{Proof}
\begin{enumerate}[(a)]	
\item
Since $(b \to c')' = b$, by Lemma \ref{pretransitivity1} (\ref{211114_07}) with $d = a'$, we have 
	$(c \to (b \to a')) \to b = (0 \to a') \to b.$
\item
	$$
	\begin{array}{lcll}
	(c \to a') \to b & = & [c \to (b \to a')] \to b & \mbox{by Lemma \ref{pretransitivity1} (\ref{031114_12})} \\
	& = & (0 \to a') \to b & \mbox{by (\ref{211114_08})}, 
	\end{array}
	$$
	from which we get
	$(c \to a') \to b = (0 \to a') \to b,$
	which, together with Lemma \ref{pretransitivity1} (\ref{031114_10}), implies 
	$(c \to a') \to b = a' \to b.$
\item
	$$
	\begin{array}{lcll}
	c \to a' & = & (0 \to a) \to (c \to a') & \mbox{by Lemma \ref{general_properties2} (\ref{281014_07}) with } x = a', y = c \\
	& = & [0 \to (a' \to b)] \to (c \to a') & \mbox{by Lemma \ref{pretransitivity1} (\ref{201114_04})} \\
	& = & [0 \to \{(c \to a') \to b\}] \to (c \to a') & \mbox{by (\ref{211114_10})} \\
	& = & [(c \to a') \to (0 \to b)] \to (c \to a') & \mbox{by Lemma \ref{general_properties2} (\ref{071114_04})} \\
	& = & [0 \to (0 \to b)] \to (c \to a') & \mbox{by Lemma \ref{general_properties2} (\ref{291014_10})} \\
	& = & (0 \to b) \to (c \to a') & \mbox{by Lemma \ref{general_properties2} (\ref{311014_03})} \\
	& = & [(c \to a') \to b] \to (c \to a') & \mbox{by Lemma \ref{general_properties2} (\ref{291014_10})} \\
	& = & (a' \to b) \to (c \to a') & \mbox{by (\ref{211114_10})}.
	\end{array}
	$$
\item
	$$
	\begin{array}{lcll}
	c \to a' & = & (0 \to a) \to (c \to a') & \mbox{by Lemma \ref{general_properties2} (\ref{281014_07})} \\
	& = & (0 \to a) \to [(a' \to b) \to (c \to a')] & \mbox{by (\ref{251114_01})} \\
	& = & (0 \to a) \to [(a' \to b) \to (c'' \to a')] & \mbox{} \\
	& = & (0 \to a) \to [\{b \to (a \to c')\}' \to (0 \to a)'] & \mbox{by Lemma \ref{transitivity} (\ref{201114_09}) with }\\
	&  &  &  d = 0, e = a, f = c'\\
	& = & (0 \to a) \to [\{b \to (a \to c')\}' \to \{0 \to (a' \to c)\}'] & \mbox{by Lemma \ref{transitivity} (\ref{201114_07})} \\
	& = & (0 \to a) \to [\{b \to (a \to c')\}' \to \{0 \to (a \to c')'\}'] & \mbox{by Lemma \ref{general_properties2} (\ref{191114_05})} \\
	& = & (0 \to a) \to [\{b \to (a \to c')\}' \to [0 \to \{b \to (a \to c')\}']'] & \mbox{by Lemma \ref{transitivity}(\ref{191114_04}) with } \\
	&  &  & d = c', e = 0 \\
	& = & [b \to (a \to c')]' \to [(a \to 0) \to \{b \to (a \to c')\}']' & \mbox{by Lemma \ref{general_properties2} (\ref{181114_16}) with }  \\
	&  &  & x = [b \to (a \to c')]', \\
	&  &  & y = a, z = 0 \\
	& = & [\{0 \to (a \to 0)\} \to \{b \to (a \to c')\}']' & \mbox{by Lemma \ref{general_properties2} (\ref{291014_09}) and (\ref{291014_10})}\\ 
	&  &  & \text{with } x = a \to 0, \\
	&  &  & y = [b \to (a \to c')]' \\
	& = & [\{a \to (0 \to 0)\} \to \{b \to (a \to c')\}']' & \mbox{by Lemma \ref{general_properties2} (\ref{071114_04})} \\
	& = & [(a \to 0') \to \{b \to (a \to c')\}']' & \mbox{} \\
	& = & [\{b \to (a \to c')\} \to [a \to \{b \to (a \to c')\}]']'  & \mbox{by Lemma \ref{transitivity}(\ref{101114_03}) with }\\
	&  &  &  e = a, d = a \to c' \\
	& = & [\{b \to (a \to c')\} \to a] \to [b \to (a \to c')] & \mbox{by Lemma \ref{general_properties2} (\ref{291014_09})} \\
	& = & a \to [b \to (a \to c')] & \mbox{by Lemma \ref{transitivity}(\ref{101114_04})}.
	\end{array}
	$$

\item
	$$
	\begin{array}{lcll}
	0 \to (a \to d) & = & 0 \to [(0 \to c) \to (a \to d)] & \mbox{by Lemma \ref{transitivity}(\ref{041114_02})} \\
	& = & c \to [0 \to (a \to d)] & \mbox{by Lemma \ref{general_properties2} (\ref{191114_02})} \\
	& = &  0 \to [c \to (a \to d)] & \mbox{by Lemma \ref{general_properties2} (\ref{071114_04})}.
	\end{array}
	$$
\item	
\noindent	
	$$
	\begin{array}{lcll}
	(d \to a) \to d & = & (0 \to a) \to d & \mbox{by Lemma \ref{general_properties2} (\ref{291014_10})} \\
	& = & [0 \to (a' \to b)] \to d & \mbox{by Lemma \ref{pretransitivity1} (\ref{201114_04})} \\
	& \sqsubseteq & (a' \to b) \to d & \mbox{by Lemma \ref{general_properties2} (\ref{011214_04})}.
	\end{array}
	$$
\item
	$$
	\begin{array}{lcll}
	(a' \to b) \to c' & = & [(c \to a') \to (b \to c')']' & \mbox{by (I)} \\
	& = & [(c \to a') \to b]' & \mbox{by hypothesis} \\
	& = & [\{c \to (b \to a')\} \to b]'  & \mbox{by Lemma \ref{pretransitivity1} (\ref{031114_12})} \\ 
	& = & [(0 \to a') \to b]' & \mbox{by Lemma \ref{pretransitivity1} (\ref{211114_07}) with } \\
	&  &  & d = a' \mbox{ since } b \sqsubseteq c \\
	& = & [(b \to a') \to b]' & \mbox{by Lemma \ref{general_properties2} (\ref{291014_10})} \\
	& = & [(b \to a') \to b'']' &  \\
	& = & [(b \to a') \to (0' \to b')']' & \mbox{by Lemma \ref{general_properties_equiv} (\ref{TXX})} \\
	& = & (a' \to 0') \to b' & \mbox{by (I)} \\
	& = & (0 \to a) \to b'. & \mbox{by Lemma \ref{general_properties} (\ref{cuasiConmutativeOfImplic2})}.
	\end{array}
	$$


	Hence, one has
	$(a' \to b) \to c' = (0 \to a) \to b'.$
\item
From Lemma \ref{transitivity} (\ref{031114_13}), we have $(0 \to c') \to b = b$.  
 Hence, we can use Lemma \ref{pretransitivity1}. 
  Therefore, we have
		
	$$
	\begin{array}{lcll}
	0 \to (a \to c') & = & 0 \to [\{(0 \to c') \to a\} \to c'] & \mbox{by Lemma \ref{pretransitivity2} (\ref{011214_10}) and Lemma \ref{transitivity}(\ref{031114_13}) } \\
	& = & [(0 \to c') \to a] \to (0 \to c') & \mbox{by Lemma \ref{general_properties2} (\ref{071114_04})} \\
	& \sqsubseteq & (a' \to b) \to (0 \to c') & \mbox{by (\ref{011214_05}) with } d = 0 \to c' \\
	& = & 0 \to [(a' \to b) \to c'] & \mbox{by Lemma \ref{general_properties2} (\ref{071114_04})} \\
	& = & 0 \to [(0 \to a) \to b'] & \mbox{by (\ref{011214_07})} \\
	& = & 0 \to [(b \to 0) \to (a \to b')']' & \mbox{by (I)} \\
	& = & 0 \to [b' \to (a \to b')']' & \mbox{} \\
	& = & 0 \to (b' \to a)' & \mbox{by hypothesis} \\
	& = & 0 \to (b \to a') & \mbox{by Lemma \ref{general_properties2} (\ref{191114_05})} \\
	& = & 0 \to a' & \mbox{by Lemma \ref{pretransitivity1} (\ref{031114_12})}.
	\end{array}
	$$
\item

\noindent 
	$$
	\begin{array}{lcll}
	0 \to a' & = & 0 \to (a \to 0) & \mbox{} \\
	& = & 0 \to [c \to (a \to 0)] & \mbox{by (\ref{011214_02})} \\
	& = & 0 \to (c \to a') & \mbox{} \\
	& = & 0 \to [(a \to c')' \to (0 \to a)'] & \mbox{by Lemma \ref{general_properties2} (\ref{271114_02})} \\
	& = & [0 \to (a \to c')]' \to (0 \to a)' & \mbox{by Lemma \ref{general_properties2} (\ref{011214_01}) and Lemma \ref{general_properties2} (\ref{311014_03})} \\
	& = & [0 \to (a \to c')]' \to (a' \to 0')' & \mbox{by Lemma \ref{general_properties} (\ref{cuasiConmutativeOfImplic2})} \\
	& = & [\{0 \to (a \to c')\} \to (0 \to a')']' & \mbox{by Lemma \ref{general_properties} (\ref{201114_08}) with } x = 0 \to (a \to c'), y = a' \\
	& = & 0 \to (a \to c') & \mbox{by (\ref{011214_08})}.
	\end{array}
	$$
\item
	
\noindent	
	$$
	\begin{array}{lcll}
	[\{c \to (a \to c')\} \to \{0 \to (a \to c')\}']' & = & [\{0 \to (a \to c')\} \to c] \to [(a \to c') \to \{0 \to (a \to c')\}']' & \\      
	&    & \hspace{4cm}             \mbox{by (I)} \\
	& = & [\{0 \to (a \to c')\} \to c] \to [\{(a \to c') \to 0\} \to (a \to c')] &  \\
	&  &   \hspace{4cm}   \mbox{by Lemma \ref{general_properties2} (\ref{291014_09})} & \mbox{} \\
	& = & [\{0 \to (a \to c')\} \to c] \to [(a \to c')' \to (a \to c')] & \mbox{} \\
	& = & [\{0 \to (a \to c')\} \to c] \to (a \to c') &  \\
	&  &   \hspace{4cm}     \mbox{by Lemma \ref{general_properties_equiv} (\ref{LeftImplicationwithtilde})} & \mbox{} \\
	& = & c \to (a \to c') & \mbox{} \\
	&  & \hspace{4cm}  \mbox{by Lemma \ref{general_properties2} (\ref{291014_08}) with } \\
	&   & \hspace{4cm}  x = a \to c', y = c.&\\
	\end{array}
	$$
\item
	
        
        
        From (\ref{021214_02}) we have that $c \to (a \to c') \sqsubseteq 0 \to (a \to c')$.  Then using (\ref{021214_01}) we get $c \to (a \to c') \sqsubseteq 0 \to a'$.
                      
\item
	$$
	\begin{array}{lcll}
	[c \to (a \to c')]' \to (0 \to a)' & = & [c \to (a \to c')]' \to (a' \to 0')' & \mbox{by Lemma \ref{general_properties} (\ref{cuasiConmutativeOfImplic2})} \\
	& = & [\{c \to (a \to c')\} \to (0 \to a')']' & \mbox{by Lemma \ref{general_properties2} (\ref{201114_08})} \\
	& = & c \to (a \to c') & \mbox{by (\ref{031214_03})}.
	\end{array}
	$$
\item
	$$
	\begin{array}{lcll}
	a \to [b \to (a \to c')] = & = & c \to a'  & \mbox{by (\ref{251114_02})} \\
	& = & c'' \to a' & \mbox{} \\
	& = & (a \to c')' \to (0 \to a)' & \mbox{by Lemma \ref{general_properties2} (\ref{271114_02})} \\
	& = & [c \to (a \to c')]' \to (0 \to a)' & \mbox{by Lemma \ref{general_properties2} (\ref{281114_01})} \\
	& = & c \to (a \to c') & \mbox{by (\ref{031214_04})} \\
	& = & a \to c' & \mbox{by Lemma \ref{general_properties2} (\ref{281114_01})}.
	\end{array}
	$$
\item
	From (\ref{251114_02}) and (\ref{031214_05}), we get
	$c \to a' = a \to c'.$
\end{enumerate}
\end{Proof}

\begin{lema} \label{transitivity1}
	Let $\mathbf A \in \mathbf I_{2,0}$ and let $a,b,c \in A$ such that $a \sqsubseteq b $ and $b \sqsubseteq c$. 
Then 
\begin{enumerate}[{\rm(a)}]	




\item 
	$c' \to [(c \to d) \to b] \sqsubseteq c$  \label{031214_08}

\item 
	$0 \to a' = c \to (0 \to a') $ \label{031214_10}

\item 
	$c' \to (a' \to b) \sqsubseteq c$ \label{031214_12}

\item 
	$(0 \to a') \to b = (c \to a') \to b$ \label{031214_18}

\item 
	$c' \to (a' \to b) \sqsubseteq 0 \to c$ \label{031214_20}

\item 
	$[(0 \to a) \to b'] \to c = c' \to (a' \to b)$ \label{031214_22}

\item
	$a' \to c = c' \to (a' \to b) $ \label{031214_23}

\item 
	$a' \to c \sqsubseteq c$ \label{031214_24}

\item 
	$a' \to c = (0 \to a') \to c$. \label{031214_25}

\end{enumerate} 	
\end{lema}

\begin{Proof}

\begin{enumerate}[(a)]
\item
	$$
	\begin{array}{lcll}
	c' \to [(c \to d) \to b] & = & c' \to [(c \to d) \to (b \to c')'] & \mbox{by hypothesis} \\
	& = & c' \to [(d \to b) \to c']' & \mbox{by (I)} \\
	& \sqsubseteq & c'' & \mbox{by Lemma \ref{general_properties2} (\ref{031214_07})} \\
	& = & c. & \mbox{}
	\end{array}
	$$
\item
	$$
	\begin{array}{lcll}
	0 \to a' & = & b \to (0 \to a') & \mbox{by Lemma \ref{pretransitivity1} (\ref{031214_11})} \\
	& = & [0 \to \{(0 \to a') \to c\}] \to [b \to (0 \to a')] & \mbox{by Lemma \ref{transitivity} (\ref{031114_13})}\\
	&  &  & \text{ and Lemma \ref{pretransitivity2} (\ref{031114_14})}\\
	&  &  &  \text{with } d = 0 \to a', e = c \\
	& = & [0 \to \{(0 \to a') \to c\}] \to (0 \to a') & \mbox{by Lemma \ref{pretransitivity1} (\ref{031214_11})} \\
	& = & [(0 \to a') \to (0 \to c)] \to (0 \to a') & \mbox{by Lemma \ref{general_properties2} (\ref{071114_04})} \\
	& = & [0 \to (0 \to c)] \to (0 \to a')  & \mbox{by Lemma \ref{general_properties2} (\ref{291014_10})} \\
	& = & (0 \to c) \to (0 \to a')  & \mbox{by Lemma \ref{general_properties2} (\ref{031114_04})} \\
	& = & c \to (0 \to a') & \mbox{by Lemma \ref{general_properties2} (\ref{311014_06})}.
	\end{array}
	$$
	
\item
	$$
	\begin{array}{lcll}
	c' \to (a' \to b) & = & c' \to [(0 \to a') \to b] & \mbox{by Lemma \ref{pretransitivity1} (\ref{031114_10})} \\
	& = & c' \to [\{c \to (0 \to a')\} \to b] & \mbox{by (\ref{031214_10})} \\
	& \sqsubseteq & c & \mbox{by (\ref{031214_08}) with } d = 0 \to a'.
	\end{array}
	$$
	
\item
	$$
	\begin{array}{lcll}
	(0 \to a') \to b & = & [c \to (0 \to a')] \to b & \mbox{by (\ref{031214_10})} \\
	& = & [(b' \to c) \to \{(0 \to a') \to b\}']' & \mbox{by (I)} \\
	& = & [(b' \to c) \to \{(b' \to 0) \to (a' \to b)'\}]' & \mbox{by (I) and $x'' \approx x$} \\
	& = & [(b' \to c) \to \{b \to (a' \to b)'\}]' & \mbox{} \\
	& = & [(b' \to c) \to \{(a' \to 0') \to b'\}]' & \mbox{by Lemma \ref{general_properties2} (\ref{031214_16})} \\
	& = & [(b' \to c) \to \{(0 \to a) \to b'\}]' & \mbox{by Lemma \ref{general_properties} (\ref{cuasiConmutativeOfImplic2})} \\
	& = & [(b' \to c) \to (a' \to b)']' & \mbox{by Lemma \ref{pretransitivity1} (\ref{031214_17})} \\
	& = & (c \to a') \to b & \mbox{by (I)}.
	\end{array}
	$$
	
\item
	$$
	\begin{array}{lcll}
	c' \to (a' \to b) & = & c' \to [(0 \to a') \to b] & \mbox{by Lemma \ref{pretransitivity1} (\ref{031214_21}) and Lemma \ref{pretransitivity1} (\ref{031214_17})} \\
	& = & c' \to [(c \to a') \to b] & \mbox{by (\ref{031214_18})} \\
	& \sqsubseteq & 0 \to c. & \mbox{by Lemma \ref{pretransitivity1} (\ref{031214_19}) with } d = a'.
	\end{array}
	$$
	
\item
	$$
	\begin{array}{lcll}
	c' \to (a' \to b) & = & [\{c' \to (a' \to b)\} \to (0 \to c)']' & \mbox{by (\ref{031214_20})} \\
	& = & [(a' \to b) \to 0] \to c & \mbox{by (I)} \\
	& = & (a' \to b)' \to c & \mbox{} \\
	& = & [(0 \to a) \to b'] \to c & \mbox{by Lemma \ref{pretransitivity1} (\ref{031214_17})}.
	\end{array}
	$$
\item
	$$
	\begin{array}{lcll}
	c' \to (a' \to b) & = & ((0 \to a) \to b') \to c & \mbox{by (\ref{031214_22})} \\
	& = & [(0 \to a) \to 0'] \to (b' \to c) & \mbox{by Lemma \ref{pretransitivity1} (\ref{031214_15}) with } d = 0 \to a \\
	& = & [(a' \to 0') \to 0'] \to (b' \to c) & \mbox{by Lemma \ref{general_properties} (\ref{cuasiConmutativeOfImplic2})} \\
	& = & [(a' \to 0) \to 0'] \to (b' \to c) & \mbox{by Lemma \ref{general_properties2} (\ref{281014_05})} \\
	& = & [a'' \to 0'] \to (b' \to c) & \mbox{} \\
	& = & (a \to 0') \to (b' \to c) & \mbox{} \\
	& = & (a \to b') \to c & \mbox{by Lemma \ref{pretransitivity1} (\ref{031214_15}) with } d = a \\
	& = & a' \to c & \mbox{by hypothesis}.
	\end{array}
	$$
\item
	This is immediate from (\ref{031214_23}) and (\ref{031214_12}). 
\item
	$$
	\begin{array}{lcll}
	(0 \to a') \to c & = & (c \to a') \to c & \mbox{by Lemma \ref{general_properties2} (\ref{291014_10})} \\
	& = & [c \to (a' \to c)']' & \mbox{by Lemma \ref{general_properties2} (\ref{291014_09})} \\
	& = & [(a' \to c) \to c']' & \mbox{by Lemma \ref{general_properties2} (\ref{071114_05})} \\
	& = & a' \to c & \mbox{by (\ref{031214_24})}.
	\end{array}
	$$

\end{enumerate}	
\end{Proof}

We are now ready to complete the proof of transitivity of $\sqsubseteq$.
\begin{Theorem} \label{order_partial_theorem}
$\sqsubseteq$ is transitive.

\end{Theorem}

\begin{Proof} Let $a,b,c \in A$ such that $a \sqsubseteq b $ and $b \sqsubseteq c$.
Observe that
	$$
	\begin{array}{lcll}
	a' & = & a \to 0 &  \\
	& = & (0 \to c) \to (a \to 0) & \mbox{by Lemma \ref{transitivity} (\ref{041114_02}) with } d = 0 \\
	& = & (0 \to c) \to a' & \mbox{} \\
	& = & (a' \to c) \to a' & \mbox{by Lemma \ref{general_properties2} (\ref{291014_10})} \\
	& = & ((0 \to a') \to c) \to a' & \mbox{by Lemma \ref{transitivity1} (\ref{031214_25})} \\
	& = & c \to a' & \mbox{by Lemma \ref{general_properties2} (\ref{291014_08})} \\
	& = & a \to c'  & \mbox{by Lemma \ref{Lemma_36} (\ref{031214_06})}.  
	\end{array}
	$$
	Consequently,
	$$a = a'' = (a \to c')',$$
	implying $a \sqsubseteq c$.  Hence, $\sqsubseteq$ is transitive on $\mathbf{A}.$
\end{Proof}

\medskip

We are now prepared to present our first main theorem.

\begin{Theorem}  \label{order_partial_theorem1}
The variety  $\mathbf{I_{2,0}}$ is a maximal subvariety of $\mathbf{I}$ with respect to the property that the relation $\sqsubseteq$ introduced in Definition \ref{order_relation} is a partial order.
\end{Theorem}

\begin{Proof}
Let $\mathbf A \in \mathbf I_{2,0}$. The relation $\sqsubseteq$ 
	is a partial order on $A$ in view of Lemma \ref{general_properties_equiv} (\ref{reflexivity}), Lemma \ref{antisymmetry}, and Theorem \ref{order_partial_theorem}.  
	
   Next, let $\mathbf V$ be a subvariety of $\mathbf I$ such that 
   $\sqsubseteq$ 
    is a partial order on every algebra in $\mathbf V$.  Now let $\mathbf A \in \mathbf V$.  Reflexivity of 
   $\sqsubseteq$ implies that $\mathbf A \models (x \to x')' \approx x$.  Therefore, by Lemma \ref{general_properties_equiv}, we conclude that $\mathbf A \in \mathbf I_{2,0}$, and hence, $ \mathbf V \subseteq \mathbf I_{2,0}$, completing the proof.   
\end{Proof}

\section{A method to construct finite $\mathbf I_{2,0}$-chains} \label{section_example_chain}

Now that we know the relation $\sqsubseteq$ is a partial order on algebras in  $\mathbf I_{2,0}$, it is natural to consider those algebras in $\mathbf I_{2,0}$, in which $\sqsubseteq$ is a total order.

\begin{definition} 
	Let $\mathbf A \in \mathbf I$. We say that $\mathbf A$ is an {\em $\mathbf I_{2,0}$-chain} (chain, for short) if  $\mathbf A \in \mathbf I_{2,0}$ and the relation $\sqsubseteq$ (see Definition \ref{order_relation})  is totally ordered on $A$.
\end{definition}
In this section we describe a method of constructing finite $\mathbf I_{2,0}$-chains. 
But, first, we will present some examples of  $\mathbf I_{2,0}$-chains that will 
foreshadow the method to construct finite $\mathbf I_{2,0}$-chains. 
 We note that, in these examples, the number $0$ is the constant element. \\

It is easy to see that the only $2$-element $\mathbf I_{2,0}$-chains, up to isomorphism,  are 

\medskip

\begin{minipage}{0.5 \textwidth}
	\begin{tabular}{r|rr}
		$\to$: & 0 & 1\\
		\hline
		0 & 1 & 1 \\
		1 & 0 & 1
	\end{tabular} with $0 \sqsubset 1$. 
	
\end{minipage}
\begin{minipage}{0.5 \textwidth}
	\begin{tabular}{r|rr}
		$\to$: & -1 & 0\\
		\hline
		-1 & -1 & -1 \\
		0 & -1 & 0
	\end{tabular} with $-1 \sqsubset 0$  
	
\end{minipage}

\bigskip

and the only $3$-element $\mathbf I_{2,0}$-chains, up to isomorphism,  are
\ \\ \ \\
\begin{minipage}{0.5 \textwidth}
	\begin{tabular}{r|rrr}
		$\to$: & 0 & 1 & 2\\
		\hline
		0 & 2 & 2 & 2 \\
		1 & 1 & 1 & 2 \\
		2 & 0 & 1 & 2
	\end{tabular} with $0 \sqsubset 1 \sqsubset 2$, 
	
\end{minipage}
\begin{minipage}{0.5 \textwidth}
	\begin{tabular}{r|rrr}
		$\to$: & -1 & 0 & 1\\
		\hline
		-1 & -1 & -1 & -1 \\
		0 & -1 & 1 & 1 \\
		1 & -1 & 0 & 1
	\end{tabular} with $-1 \sqsubset 0 \sqsubset 1$, 
	
\end{minipage}

\medskip

\begin{minipage}{0.5 \textwidth}
	\begin{tabular}{r|rrr}
		$\to$: & -2 & -1 & 0\\
		\hline
		-2 & -2 & -2 & -2 \\
		-1 & -2 & -1 & -1 \\
		0 & -2 & -1 & 0
	\end{tabular} with $-2 \sqsubset -1 \sqsubset 0$. 
	
\end{minipage}
\medskip




\medskip 

Note that, henceforth, we will use the symbol $\leq$ to denote the natural order in $\mathbb Z$. Recall that $\sqsubseteq$ is being used for the order (see Definition \ref{order_relation}).

\medskip

The next definition describes a general method to construct a class of finite $\mathbf{I_{2,0}}$-chains, generalizing the above examples.  In the next section, we will show that, this method, in fact, yields, up to isomorphism, all finite $\mathbf{I_{2,0}}$-chains.  

\begin{definition}\label{definition_functions_pandast}
	Let $k \in \mathbb{N}$.   Let $m, n \in \omega$ be such that the interval $[-n, m] \subseteq \mathbb Z$ with $|[-n, m]| = k$ and $0 \leq n, m \leq k-1$. The (auxiiliary) functions $p$ (predecessor) and $\ast$ are defined on $[-n, m]$ as follows:
	$$p(x) = \left\{\begin{array}{lcl}
	x-1 & \mbox{if} & x > -n \\
	-n & \mbox{if} & x = -n,
	\end{array}
	\right.$$
	and
	$$x^\ast = \left\{\begin{array}{lcl}
	m  & \mbox{if} & x = 0 \\
	x  & \mbox{if} & x<0 \\
p((p(x))^\ast)  & \mbox{if} & x > 0.
	\end{array}
	\right.$$
For convenience, we write $p(p(x)^\ast)$ for $p((p(x))^\ast)$.	\rm{(}Notice that the function $\ast$ is defined recursively for $x \geq 0$.\rm{)}\\
Define the algebra  $\mathbf{[-n, m]}$ as follows:
\begin{center}
 $\mathbf{[-n, m]} := \langle [-n, m]; \Rightarrow, 0 \rangle$, where $0 \in [-n, m]$ is the constant and $\Rightarrow$ is defined by
	$$x \Rightarrow y = \left\{\begin{array}{ll}
	max(x^\ast,y)  & \mbox{if } x,y \geq 0   \\
	min(x,y)  & \mbox{otherwise.}
	\end{array}
	\right.$$
\end{center}
	We set $x' := x \Rightarrow  0$. 

\end{definition}


\medskip

We shall now illustrate the method described in the above definition by applying it to construct a $6$-element $\mathbf{I_{2,0}}$-chain.\\

Let $k = 6$, and consider the interval $ A =  [-2, 3] = \{-2,-1,0,1,2,3\}$. Since $0 \Rightarrow 0 = max(0^\ast,0) = max(3,0) = 3$ and $a \Rightarrow b = min(a,b)$ if $a<0$ or $b < 0$, we arrive at tha following partial table for $\Rightarrow$:\\

\medskip

\noindent \begin{tabular}{c |c c c c c c }
	\hline
	$\Rightarrow$ &  -2 &  -1 &  0 &  1 &  2 & 3  \\ \hline
	-2 &  -2 &  -2 &  -2 &  -2 &  -2 &  -2 \\
	-1 &  -2 &  -1 &  -1 &  -1 &  -1 & -1  \\
	0 &  -2 &  -1 &  3 &  ? &  ? & ?  \\
	1 &  -2 &  -1 &  ? &  ? &  ? & ?  \\
	2 &  -2 &  -1 &  ? &  ? &  ? & ?  \\
	3 &  -2 &  -1 &  ? &  ? &  ? & ?  \\
\end{tabular} \\

\medskip

Next, we determine the operations $p$ and $\ast$:

\medskip

\noindent \begin{minipage}{0.5 \textwidth}
	\noindent \begin{tabular}{c |l }
		\hline
		$x$ &  $x^\ast$ \\ \hline
		$0$ & $3$  \\ 
		$1$ & $p(p(1)^\ast) = p(0^\ast) = p(3) = 2$  \\ 
		$2$ & $p(p(2)^\ast) = p(1^\ast) = p(2) = 1$  \\ 
		$3$ & $p(p(3)^\ast) = p(2^\ast) = p(1) = 0$  
	\end{tabular}
\end{minipage}
\begin{minipage}{0.5 \textwidth}
	\noindent \begin{tabular}{c |l }
		\hline
		$x$ &  $x \Rightarrow 0$ \\ \hline
		$1$ & $max(1^\ast,0) = max(2,0) = 2$  \\ 
		$2$ & $max(2^\ast,0) = max(1,0) = 1$  \\ 
		$3$ & $max(3^\ast,0) = max(0,0) = 0$  
	\end{tabular}
\end{minipage}

%


\medskip

Hence the table for $\Rightarrow$ becomes:

\medskip

\noindent \begin{tabular}{c |c c c c c c }
	\hline
	$\Rightarrow$ &  -2 &  -1 &  0 &  1 &  2 & 3  \\ \hline
	-2 &  -2 &  -2 &  -2 &  -2 &  -2 &  -2 \\
	-1 &  -2 &  -1 &  -1 &  -1 &  -1 & -1  \\
	0 &  -2 &  -1 &  3 &  ? &  ? & ?  \\
	1 &  -2 &  -1 &  2 &  ? &  ? & ?  \\
	2 &  -2 &  -1 &  1 &  ? &  ? & ?  \\
	3 &  -2 &  -1 &  0 &  ? &  ? & ?  \\
\end{tabular}

\medskip

Observe that $0 \Rightarrow 1 = max(0^\ast,1) = max(3,1) = 3$, $1 \Rightarrow 1 = max(1^\ast,1) = max(2,1) = 2$, $2 \Rightarrow 1 = max(2^\ast,1) = max(1,1) = 1$ and $3 \Rightarrow 1 = max(3^\ast,1) = max(0,1) = 1$. Then
we get\\
\medskip

\noindent \begin{tabular}{|c |c |c |c |c |c |c |}
	\hline
	$\Rightarrow$ &  -2 &  -1 &  0 &  1 &  2 & 3  \\
	-2 &  -2 &  -2 &  -2 &  -2 &  -2 &  -2 \\
	-1 &  -2 &  -1 &  -1 &  -1 &  -1 & -1  \\
	0 &  -2 &  -1 &  3 &  3 &  ? & ?  \\
	1 &  -2 &  -1 &  2 &  2 &  ? & ?  \\
	2 &  -2 &  -1 &  1 &  1 &  ? & ?  \\
	3 &  -2 &  -1 &  0 &  1 &  ? & ?  \\
	\hline
\end{tabular}\\

\medskip

Iterating this process we obtain the following complete table for $\Rightarrow$:\\

\medskip

\noindent \begin{tabular}{|c |c |c |c |c |c |c |}
	\hline
	$\Rightarrow$ &  -2 &  -1 &  0 &  1 &  2 & 3  \\
	-2 &  -2 &  -2 &  -2 &  -2 &  -2 &  -2 \\
	-1 &  -2 &  -1 &  -1 &  -1 &  -1 & -1  \\
	0 &  -2 &  -1 &  3 &  3 &  3 & 3  \\
	1 &  -2 &  -1 &  2 &  2 &  2 & 3  \\
	2 &  -2 &  -1 &  1 &  1 &  2 & 3  \\
	3 &  -2 &  -1 &  0 &  1 &  2 & 3  \\
	\hline
\end{tabular}\\

Thus we have constructed the algebra $\mathbf{[-n, m]}$.  Observe that $-2 \sqsubset -1 \sqsubset 0  \sqsubset 1 \sqsubset 2 \sqsubset 3$ and $x'' = x^{**} = x$.  Also, it is routine to verify $\mathbf{[-n, m]} \in \mathbf{I_{2,0}}$. Hence it is an $\mathbf I_{2,0}$-chain. 

\medskip
\ \par
Returning to the general method, we now aim to prove that $\mathbf{[-n;m]}$
is an $\mathbf I_{2,0}$-chain. To prove this, we will need the following lemmas.

\begin{lema} \label{161214_01}
	If $x \in \mathbf{[-n, m]}$ and $0 \leq x \leq m$ then $x^\ast = m - x$ and, consequently, $x^\ast \in [0, m]$.
\end{lema}

\begin{Proof}
	We prove this lemma by induction on the element $x$. Assume that $x = 0$. Then $0^\ast = m = m-0$. 
	
	Next, suppose $x > 0$. Since $-n \leq 0 < x$, we have $p(x) = x-1$.
	Hence, by inductive hypothesis, we have 
	\begin{equation}\label{171214_01}
	p(x)^\ast = m - p(x) = m - (x-1) = m-x+1.
	\end{equation}
	From $x > 0$, we can conclude that $m-x+1 \leq m$. Also, since $x \leq m$, we obtain $0 \leq m-x$, thus  $-n-1 < 0 \leq m-x$, implying $m-x+1 > -n$. Then we get $p(m-x+1) = m-x+1-1$. 
	By (\ref{171214_01}), $x^\ast = p((p(x))^\ast) = p(m-x+1) = m-x$, completing the induction.  It is clear that $x^\ast \in [0, m]$.  
\end{Proof}

\begin{corolario} \label{171214_02}
	If $x \in \mathbf{[-n, m]}$ then $x' = x^\ast$.
\end{corolario}

\begin{Proof}
If $x < 0$ we have that $x' = x \Rightarrow 0 = min(x,0) = x = x^\ast$.  If $x > 0$, then by Lemma \ref{161214_01}, $x^\ast \geq 0$, and hence $x' = x \Rightarrow 0 = max(x^\ast,0) = x^\ast$.
\end{Proof}

\begin{lema} \label{161214_02}
	If $x \in \mathbf{[-n, m]}$ then $x'' = x$.
\end{lema}

\begin{Proof}
We consider the following cases:
	\begin{itemize}
		\item If $x < 0$, then $x^\ast =x$, and hence  $x^{\ast \ast}= x$.
		\item If $x \geq 0$,
		$$
		\begin{array}{lcll}
		x^{\ast \ast} & = & (m-x)^\ast & \mbox{by Lemma \ref{161214_01} since } 0 < x \leq m \\
		& = & m - (m-x) & \mbox{by Lemma \ref{161214_01} since } 0 \leq m-x \leq m \\
		& = & x. & \mbox{}
		\end{array}
		$$
	\end{itemize}
	Consequently, by Corollary \ref{171214_02}, $x'' = x$.
\end{Proof}

\begin{lema} \label{171214_03}
	If $x,y \in \mathbf{[-n, m]}$ and $0 \leq x \leq y$ then $x^\ast \geq y^\ast$.
\end{lema}

\begin{Proof}
	We prove this lemma by induction on the element $x$. If $x = 0$, $x^\ast = 0^\ast = m \geq y^\ast$ by Lemma \ref{161214_01}.
	
	Now assume that $x > 0$. Since $0 < x \leq y$, we have that $x^\ast = p(p(x)^\ast)$ and $y^\ast = p(p(y)^\ast)$. Note that $0 \leq p(x) \leq p(y)$. Then, by induction hypothesis, we get $p(y)^\ast \leq p(x)^\ast$. Hence $x^\ast = p(p(x)^\ast) \geq p(p(y)^\ast) = y^\ast$.
\end{Proof}

\medskip


\begin{lema}\label{example_I20_chain}
Let $k \in \mathbb{N}$.   Let $m, n \in \omega$ be such that the interval $[-n, m] \subseteq \mathbb Z$ with $|[-n, m]| = k$ and $0 \leq n, m \leq k-1$.  Then, $\mathbf{ [-n, m]} \in \mathbf I_{2,0}$.   
\end{lema}


\begin{Proof} 
	The proof that $\langle [-n;m]; \Rightarrow, 0 \rangle$ satisfies the identity (I) is long and computational, but routine.  Hence we leave the verification to the reader with the recommendation that the following cases be considered, where $i,j,k \in [-n;m]$:
	
	\begin{multicols}{2}
		\begin{enumerate}[(1)]
			\item $i,j,k \geq 0$, $i^* \geq j$, $i \geq k$  
			\item $i,j,k \geq 0$, $i^* \geq j$, $i < k$    
			\item $i,j,k \geq 0$, $i^* < j$, $k \geq i$    
			\item $i,j,k \geq 0$, $i^* < j$, $k < i$, $j^* \leq k$  
			\item $i,j,k \geq 0$, $i^* < j$, $k < i$, $j^* > k$   
			\item $i,j \geq 0$ and $k<0$   
			\item $i \geq 0$,  $j<0$ and $k \geq 0$  
			\item $i \geq 0$,  $j<0$ and $k < 0$  
			\item $i<0$,  $j \geq 0$ and $k \geq 0$  
			\item $i<0$, $j \geq 0$ and $k < 0$  
			\item $i<0$, $j < 0$ and $k \geq 0$  
			\item $i,j,k<0$. 
		\end{enumerate}
	\end{multicols}
	Observe that, if $x \in [-n, m]$, then, from Corollary \ref{171214_02}, we have $x' = x^*$, and from Lemma \ref{161214_02} we have that $x'' = x$; and in particular $0''=0$.  Thus, we conclude that $\langle [-n, m]; \Rightarrow, 0 \rangle \in \mathbf I_{2,0}$. 
\end{Proof}

\medskip	
	In view of the above lemma and Theorem  \ref{order_partial_theorem}, the relation defined by 
	     $$x \sqsubseteq y \quad \mbox{if and only if} \quad (x \Rightarrow y')'=x$$
	     is a partial order on $\mathbf{ [-n, m]}$.
We now wish to show that $\sqsubseteq$ is indeed a total order.
\begin{Lemma}	\label{Lemma_order_isomorphism}
Let $\mathbf{[-n, m]}$ be the algebra, as defined in Definition \ref{definition_functions_pandast}.  Then 
$$\langle [-n, m]; \sqsubseteq \rangle \cong \langle [-n, m]; \leq \rangle.$$
\end{Lemma}	
	\begin{Proof}
	Let $x, y \in [-n, m]$.
	 It is enough to prove that $x \leq y$ if and only if $x \sqsubseteq y$.
	
	Assume that $x \leq y$. We will consider the following cases:
	\begin{itemize}
		\item {\bf Case 1}: $x < 0$. Then
		\begin{equation}\label{191214_01}
		(x \Rightarrow y')' = (x \Rightarrow y^\ast)^\ast = [min(x,y^\ast)]^\ast.
		\end{equation}
		We consider further the following subcases:
	
		\begin{itemize}
		\item {\bf Case 1.1}: $y < 0$.
		$$
		\begin{array}{lcll}
		(x \Rightarrow y')' & = & [min(x,y^\ast)]^\ast & \mbox{by (\ref{191214_01})} \\
		& = & [min(x,y)]^\ast & \mbox{since $y < 0$} \\
		& = & x^\ast & \mbox{since $x \leq y$} \\
		& = & x. & \mbox{since $x < 0$}
		\end{array}
		$$
		
		\item {\bf Case 1.2}: $y \geq 0$.
		$$
		\begin{array}{lcll}
		(x \Rightarrow y')' & = & [min(x,y^\ast)]^\ast & \mbox{by (\ref{191214_01})} \\
		& = & x^\ast & \mbox{since 
		                            $y^\ast \geq 0$ by Lemma \ref{161214_01}, and $x < 0$} \\
		& = & x. & \mbox{}
		\end{array}
		$$
		\end{itemize}
		\item {\bf Case 2}: $x \geq 0$. Therefore $y \geq 0$. In this case
		$$
		\begin{array}{lcll}
		(x \Rightarrow y')' & = &  (x \Rightarrow y^\ast)^\ast & \mbox{} \\
		& = & [max(x^\ast,y^\ast)]^\ast & \mbox{} \\
		& = & x^{\ast \ast} & \mbox{by Lemma \ref{171214_03}} \\
		& = & x & \mbox{}
		\end{array}
		$$
			
	\end{itemize}
	Thus, in all these cases, $x \sqsubseteq y$. \\
	
	For the converse, suppose $x \sqsubseteq y$.
	
	\begin{itemize}
		\item {\bf Case 1}: $x < 0$.  If $y \geq 0$ then $x < y$.  So, we can assume $y < 0$. Then
		$$
		\begin{array}{lcll}
		x & = & x' & \mbox{since $x < 0$}\\
		& = & (x \Rightarrow y')'' & \mbox{by hypothesis} \\
		& = & x \Rightarrow y' & \mbox{by Lemma \ref{161214_02}} \\
		& = & x \Rightarrow y & \mbox{} \\
		& = & min(x,y). & \mbox{}
		\end{array}
		$$
		Hence $x \leq y$.
		
		\item {\bf Case 2}: $x \geq 0$. Suppose $y<0$.  Then
		$$
		\begin{array}{lcll}
		x & = & (x \Rightarrow y')' & \mbox{by hypothesis} \\
		& = & (x \Rightarrow y)' & \mbox{} \\
		& = & min(x,y)' & \mbox{} \\
		& = & y' & \mbox{} \\
		& = & y, & \mbox{}
		\end{array}
		$$
		a contradiction. Hence $y \geq 0$. Consequently,
		$$
		\begin{array}{lcll}
		x' & = & (x \Rightarrow y')'' &  \\
		& = & x \Rightarrow y' & \mbox{by Lemma \ref{161214_02}} \\
		& = & max(x',y'), & \mbox{}
		\end{array}
		$$
	so, $x' \geq y'$. Then, by Lemma \ref{161214_02} and Lemma \ref{171214_03}, $x = x'' \leq y'' = y$.
	\end{itemize}
\end{Proof}	

  In view of Lemma \ref{example_I20_chain} and Lemma \ref{Lemma_order_isomorphism}, we have proved the following
\begin{Theorem} 
$\mathbf{[-n, m]}$ is an $\mathbf{I_{2,0}}$-chain, where
 $$-n \sqsubset -n+1 \sqsubset \ldots \sqsubset -1 \sqsubset 0 \sqsubset 1 \sqsubset 2 \sqsubset \ldots \sqsubset m.$$ 
\end{Theorem}

\section{Characterization of finite $\mathbf I_{2,0}$-chains} \label{section_chains}

\medskip

In this section we are going to prove our second main result.  
The following lemmas will be useful later in this section. 

\begin{lema} \label{last_element}
	Let $\mathbf A \in \mathbf I_{2,0}$. 
	Then $0'$ is the greatest element in $A$, relative to $\sqsubseteq$.
\end{lema}

\begin{Proof}
Let $a \in A$.  Since $(a \to (0 \to 0)')' = (a \to 0'')' = (a \to 0)' = a'' = a$, we have $a \sqsubseteq 0'$.
\end{Proof}

\begin{lema} \label{121214_09}
	Let $\mathbf A \in \mathbf I_{2,0}$ and let $a,b \in A$ with $0 \sqsubseteq a \sqsubseteq b$. Then $b' \sqsubseteq a'$.
\end{lema}

\begin{Proof}
	$$
	\begin{array}{lcll}
	(b' \to a'')'  & = & (b' \to a)' & \mbox{} \\
	& = & (b' \to (a \to b')')' & \mbox{by hypothesis}\\          
	& = & ((a \to 0') \to b'')' & \mbox{by Lemma \ref{general_properties2} (\ref{031214_16})} \\
	& = & ((a \to 0') \to b)' & \mbox{} \\
	& = & ((0 \to a') \to b)' & \mbox{by Lemma \ref{general_properties} (\ref{cuasiConmutativeOfImplic2})} \\
	& = & ((0 \to a')'' \to b)' & \mbox{} \\
	& = & (0' \to b)' & \mbox{by hypothesis}\\                                    
	& = & b' & \mbox{by Lemma \ref{general_properties_equiv} (\ref{TXX}).}
	\end{array}
	$$
	
\end{Proof}

\begin{lema} \label{121214_01}
	Let $\mathbf A \in \mathbf I_{2,0}$ and let $a \in A$. If  $0 \sqsubseteq a$ then $0 \to a = 0'$.
\end{lema}

\begin{Proof}
	First notice that, since $0 \sqsubseteq a$, $0' = (0 \to a')'' = 0 \to a'$. Consequently,
	\begin{equation} \label{121214_07}
	0' = 0 \to a'.
	\end{equation}
	Then
	$$
	\begin{array}{lcll}
	0' & = & 0' \to 0' & \mbox{by Lemma \ref{general_properties_equiv} (\ref{TXX})} \\
	& = & (0 \to a') \to 0' & \mbox{by (\ref{121214_07})} \\
	& = & (0' \to a') \to 0' & \mbox{by Lemma \ref{general_properties2} (\ref{291014_10})} \\
	& = & a' \to 0' & \mbox{by Lemma \ref{general_properties_equiv} (\ref{TXX})} \\
	& = & 0 \to a. & \mbox{by Lemma \ref{general_properties} (\ref{cuasiConmutativeOfImplic2})}
	\end{array}
	$$
\end{Proof}

\begin{lema} \label{121214_02}
	Let $\mathbf A \in \mathbf I_{2,0}$ and let $a,b \in A$. If  $0 \sqsubseteq a$ and $0 \sqsubseteq b$ then $0 \sqsubseteq a \to b $.
\end{lema}

\begin{Proof}
	$$
	\begin{array}{lcll}
	[0 \to (a \to b)']'  & = & [(a \to b) \to 0']' & \mbox{by Lemma \ref{general_properties} (\ref{cuasiConmutativeOfImplic2})} \\
	& = & (0 \to a) \to (b \to 0')' & \mbox{by (I)} \\
	& = & (0 \to a) \to (0 \to b')' & \mbox{by Lemma \ref{general_properties} (\ref{cuasiConmutativeOfImplic2})} \\
	& = & (0 \to a) \to 0 & \mbox{since }  0 \sqsubseteq  b \\
	& = & 0' \to 0 & \mbox{by Lemma \ref{121214_01} since } 0 \sqsubseteq  a \\
	& = & 0. & \mbox{by Lemma \ref{general_properties_equiv} (\ref{TXX})}
	\end{array}
	$$
\end{Proof}

\begin{corolario} \label{121214_08}
	Let $\mathbf A \in \mathbf I_{2,0}$ and $a \in A$. If  $a \sqsupseteq 0$  then $a' \sqsupseteq 0$.
\end{corolario}

\begin{lema} \label{tildeInverImplic}
	Let $\mathbf A$ be an $\mathbf I_{2,0}$-chain and let $a,b \in A$. Then $a' \to b' = b \to a$.
\end{lema}

\begin{Proof}
	Since $\mathbf A$ is a chain, we can assume that $b' \sqsubseteq a$ or $a \sqsubseteq b'$.
	
	If $b' \sqsubseteq a$, $(b' \to a')' = b'$, then $b' \to a' = b$. Hence $b \to a = (b' \to a') \to a = [(a' \to b') \to (a' \to a)']'$, using (I).  By Lemma \ref{general_properties_equiv} (\ref{LeftImplicationwithtilde}), $[(a' \to b') \to (a' \to a)']' = [(a' \to b') \to a']' = [[(a \to a') \to (b' \to a')']']' = (a \to a') \to (b' \to a')' = (a'' \to a') \to (b' \to a')' = a' \to b'$.
	
If $a \sqsubseteq b'$ then we have $a' = (a \to b'')'' = a \to b$, and the rest of the argument is similar to the previous case. 
\end{Proof}

\begin{lema} \label{lemma_properties_chains}
	Let $\mathbf A$ be a $\mathbf I_{2,0}$-chain with $|A| \geq 2$ and let $a \in A$ such that $a\sqsubset 0$. Then    
	\begin{enumerate}[{\rm (a)}]
		\item $0 \to a' = a'$ \label{271014_05}
		\item $0 \to a = a$ \label{271014_01}
		\item $(a \to a) \to a = a \to a$ \label{271014_06}
		\item $a \to a = a'$ \label{271014_07}
		\item $a \to a = a$ \label{271014_09}
		\item $a = a'$. \label{271014_10}
	\end{enumerate}
\end{lema}

\begin{Proof}
	\begin{itemize}
		\item[(\ref{271014_05})]
		Since $a \sqsubseteq 0$, we have that $a = (a \to 0')'$. Therefore, $a' = (a \to 0')'' = a \to 0' = 0 \to a'$ by Lemma \ref{general_properties}  (\ref{cuasiConmutativeOfImplic}).

		\item[(\ref{271014_01})]
		Since $a \sqsubseteq 0$, we have 
		\begin{equation} \label{271014_02}
		a = (a \to 0')'.
		\end{equation}
		Then we get
		$$\begin{array}{lcll}
		(0 \to a) \to 0'  & =  & [(0 \to 0) \to (a \to 0')']'  &  \mbox{by (I)} \\
		& =  & [(0 \to 0) \to a]'  & \mbox{by (\ref{271014_02})}  \\
		& =  & [0' \to a]'  & \mbox{}\\
		& =  & a'  & \mbox{by lemma \ref{general_properties_equiv} (\ref{TXX})}  \\
		\end{array}
		$$
		Using Lemma \ref{general_properties}  (\ref{cuasiConmutativeOfImplic}), we obtain
		\begin{equation} \label{271014_03}
		a' = 0 \to (0 \to a)'.
		\end{equation}
		Since $\mathbf A$ is a chain, $0 \sqsubseteq 0 \to a$ or $0 \to a \sqsubseteq 0$.  Suppose that $0 \sqsubseteq 0 \to a$. Then $(0 \to (0 \to a)')' = 0$ Therefore, by (\ref{271014_03}), 
		$a = a'' = (0 \to (0 \to a)')' = 0$, a contradiction, since $a \not= 0$.  Consequently, $0 \to a \sqsubseteq 0$. 
		Hence, we have
		$$\begin{array}{lcll}
		0 \to a  & =  &  ((0 \to a) \to 0')' &  \mbox{since $0 \to a \sqsubseteq 0$}\\                
		& =  & (0 \to (0 \to a)')'  & \mbox{by lemma \ref{general_properties}  (\ref{cuasiConmutativeOfImplic})}  \\
		& =  & a''  &  \mbox{by (\ref{271014_03})} \\
		& =  & a.  &
		\end{array}
		$$
		
		\item[(\ref{271014_06})]
		$$\begin{array}{lcll}
		a \to a  & =  & (0 \to a) \to a  &  \mbox{by item (\ref{271014_01})} \\
		& =  & (a' \to 0') \to a  &  \mbox{by lemma \ref{general_properties}  (\ref{cuasiConmutativeOfImplic2})} \\
		& =  & [(a' \to a') \to (0' \to a)']'  &  \mbox{by (I)} \\
		& =  & [(a \to a) \to (0' \to a)']'  &  \mbox{by Lemma \ref{tildeInverImplic}} \\
		& =  & [(a \to a) \to a']'  &  \mbox{by lemma \ref{general_properties}  (\ref{TXX})}  \\
		& =  & [[(a'' \to a) \to (a \to a')']']'  &  \mbox{by (I)} \\
		& =  & (a \to a) \to (a \to a')'  &   \\
		& =  & (a \to a) \to (a'' \to a')'  &   \\
		& =  & (a \to a) \to a''  &  \mbox{by lemma \ref{general_properties}  (\ref{LeftImplicationwithtilde})} \\
		& =  & (a \to a) \to a.  &
		\end{array}
		$$
		
		\item[(\ref{271014_07})]  
		Since $\mathbf A$ is a chain, $0 \to a'  \sqsubseteq a$ or $a \sqsubseteq 0 \to a' $.       
		
		First, we assume that $0 \to a' \sqsubseteq a$. Then
		$$\begin{array}{lcll}
		a \to a  & =  & (a \to a) \to a  &  \mbox{by (\ref{271014_06})} \\
		& =  & (a' \to a') \to a  &  \mbox{by Lemma \ref{tildeInverImplic}} \\
		& =  &  a' \to (a' \to a')' &  \mbox{by Lemma \ref{tildeInverImplic}} \\
		& =  &  (a \to 0) \to (a' \to a')' &  \mbox{} \\
		& =  & [(a \to 0) \to (a' \to a')']''  &  \mbox{}  \\
		& =  & [(0 \to a') \to a']'  &  \mbox{using (I)} \\
		& =  & 0 \to a'  &  \mbox{since $0 \to a' \sqsubseteq a$}\\                   
		& =  & a'  &  \mbox{using (\ref{271014_05})}.
		\end{array}
		$$
		Next, we assume $a \sqsubseteq 0 \to a' $, i.e.,
		$(a \to (0 \to a')')' = a$.
		Then, from (\ref{271014_05}), \\
		we have $a \to a 
		= (a \to a'')'' = [(a \to (0 \to a')')']' = a'$.

		\item[(\ref{271014_09})] Using the items (\ref{271014_06}), (\ref{271014_07}) and Lemma \ref{general_properties}  (\ref{LeftImplicationwithtilde}), we have 
		$a \to a = (a \to a) \to a = a' \to a = a$.
		
		\item[(\ref{271014_10})] This follows immediately from the two preceding items.	
	\end{itemize}
\end{Proof}

\begin{lema}
	Let $\mathbf A$ be an $\mathbf I_{2,0}$-chain with  $|A| \geq 2$, and let $a,b \in A$. If $0 \sqsubseteq a$ and $b \sqsubset 0$ then $b \to a = b$ and $a \to b = b$.
\end{lema}

\begin{Proof}
	Since $0 \sqsubseteq a$ and $b \sqsubset 0$, we have that $(0 \to a')'=0$ and $(b \to 0')' = b$. Therefore, using Lemma \ref{tildeInverImplic}, $b = b'' = b' \to 0 = (b \to 0')'' \to (0 \to a')' = (b \to 0') \to (0 \to a')' = (0 \to b') \to (a \to 0')' = [(b' \to a) \to 0']'$. Hence,
	\begin{equation} \label{271014_11}
	b = [(b' \to a) \to 0']'.
	\end{equation}
	From the hypothesis and Lemma \ref{lemma_properties_chains} (\ref{271014_10}), we have
	\begin{equation} \label{271014_13}
	b' = b.
	\end{equation}
	Suppose that $0 \sqsubseteq b' \to a$. Then $0 = [0 \to (b' \to a)']' = [(b' \to a) \to 0']'$ by Lemma \ref{tildeInverImplic}, implying  
	$0 = b$, which is a contradiction in view of (\ref{271014_11}).  Consequently, $b' \to a \sqsubseteq 0$, since $\mathbf A$ is a chain.  Hence,
	\begin{equation} \label{271014_12}
	b' \to a = [(b' \to a) \to 0']'.
	\end{equation}
	From (\ref{271014_11}), (\ref{271014_13}) and (\ref{271014_12}) we conclude 
	$b = b \to a,$ proving the first half of the conclusion of the lemma.
	From 
	$$\begin{array}{lcll}
	b  & =  & (b \to a')'  &  \mbox{since } b \sqsubseteq a, \mbox{ as } 0 \sqsubseteq a \mbox{ and }  b \sqsubset 0 \\
	& =  & (a'' \to b')'  &  \mbox{by Lemma \ref{tildeInverImplic} } \\
	& =  & (a \to b')'  &  \mbox{} \\
	& =  & (a \to b)'  &  \mbox{by } (\ref{271014_13})
	\end{array}
	$$
	we conclude that $a \to b = b' = b$ in view of (\ref{271014_13}), completing the second half.		
\end{Proof}


\begin{definition}\label{Definition_finite_algebra_as_interval}
	Let $\mathbf A = \langle A; \to, 0 \rangle$ be a finite $\mathbf I_{2,0}$-chain. We let $A^+ := \{a \in A:\ a \sqsupset 0\}$ and $A^- := \{a \in A:\ a \sqsubset 0\}$. Observe that $A = A^+ \cup \{0\} \cup A^-$. 
	Henceforth, without loss of generality, we will represent $A = [-n, m]$ with  $0 \leq n, m \leq |A|-1$, 
	 such that
	$$-n \sqsubset -n+1 \sqsubset \ldots \sqsubset -1 \sqsubset 0 \sqsubset 1 \sqsubset 2 \sqsubset \ldots \sqsubset m.$$
\end{definition}

\begin{remark}
In view of the above definition, we can use the functions $\ast$ and $p$ of Definition \ref{definition_functions_pandast} as functions on the domain $[-n,m]$ of $\mathbf{A}$ as well.  
\end{remark}
Now, we wish to prove that $\langle A; \to, 0 \rangle = \langle [-n;m]; \Rightarrow, 0 \rangle$.  To achieve this, 
we need the following lemmas.

\begin{lema} \label{lemma_recurs_p}
	Let $\mathbf A = \langle A; \to, 0 \rangle$ be a finite $\mathbf I_{2,0}$-chain with  $|A| \geq 2$. If $a \sqsupset 0$ then $a' = p(p(a)')$.
\end{lema}

\begin{Proof}
	By hypothesis we have that $a \sqsupset 0$. Then $p(a) \sqsupseteq 0$. Hence $0 \sqsubseteq p(a) \sqsubset a$. Then, by Lemma \ref{121214_09},
	\begin{equation} \label{121214_10}
	a' \sqsubseteq p(a)'.
	\end{equation}
	Since $a \sqsupset 0$, by Corollary \ref{121214_08}, $a' \sqsupseteq 0$. Therefore, by (\ref{121214_10}),
	\begin{equation} \label{121214_11}
	0 \sqsubseteq p(a)'.
	\end{equation}
	If $a' = p(a)'$ then $a = p(a)$ and, consequently, $a = -n$, a contradiction, so $a' \sqsubset p(a)'$, and hence, $0 \sqsubseteq a' \sqsubseteq p(p(a)') \sqsubset p(a)'$. By lemma \ref{121214_09}, $a \sqsupseteq [p(p(a)')]' \sqsupseteq p(a)$.   Thus
	\begin{equation} \label{121214_12}
	[p(p(a)')]' \in \{a, p(a)\}.
	\end{equation}
	If $[p(p(a)')]' = p(a)$, we have that $p(p(a)') = [p(p(a)')]'' = p(a)'$, a contradiction, since $p(a)' \sqsupseteq 0$ by (\ref{121214_11}).
	Therefore $[p(p(a)')]' = a$ and therefore, $p(p(a)') = a'$.	
\end{Proof}

\begin{lema} \label{lemma_ast_tilde}
	Let $\mathbf A = \langle A; \to, 0 \rangle$ be a finite $\mathbf I_{2,0}$-chain. If $a \in A$ then $a^\ast = a'$.
\end{lema}

\begin{Proof}
	The statement $0' = m = 0^\ast$ follows from Lemma \ref{last_element}.  If $a \sqsubset 0$ then $a' = a$ 
by Lemma \ref{lemma_properties_chains} (\ref{271014_10}), and $a= a^\ast$ by definition, implying $a= a^\ast$.  
	
	Now assume that $a \sqsupset 0$. We will verify that $a' = a^\ast$ by induction on $a$. If $a = 1$, then, as $0' = 0^\ast$, we have, by Lemma \ref{lemma_recurs_p}, that  $1' = p(p(1)') = p(0') = p(0^\ast) = p(p(1)^\ast) = 1^\ast$. The inductive hypothesis is that $p(a)' = p(a)^\ast$. Hence, we have, by Lemma \ref{lemma_recurs_p}, $a' = p(p(a)') = p(p(a)^\ast) = a^\ast$.
\end{Proof}
\medskip

The following theorem shows that the general method described in Definition 
\ref{definition_functions_pandast}     essentially gives all finite $\mathbf I_{2,0}$-chains. 

\begin{teorema} \label{theorem_chains_implic}
	Let $\mathbf A$ 
be a finite $\mathbf I_{2,0}$-chain. Then $\mathbf A \cong \langle [-n, m]; \Rightarrow, 0 \rangle$  for some $0 \leq n, m \leq |A|-1$.  
\end{teorema}


\begin{Proof} We will use the notation of Definition \ref{Definition_finite_algebra_as_interval}.  
Let $i,j \in A$. From Lemma \ref{lemma_ast_tilde}, $i' = i^\ast$ and $j' = j^\ast$.  It suffices to verify that 
	$$i \to j = \left\{\begin{array}{lc}
	 max (i' , j) & \mbox{if }  i,j \sqsupseteq 0 \\
	 min (i, j) & \mbox{otherwise}
	\end{array}
	\right.$$ with $0' = m$.  We consider the following cases:

\begin{itemize}
		\item {\bf Case 1}: $j > 0$. \\
		We need the following subcases:
		\begin{itemize}		
		\item {\bf Case 1.1}: $i > 0$. \\
		 We make the following further subcases: 
		\begin{itemize}
		\item {\bf Case 1.1.1}: $i' \geq j$.\\
		Since $i' \sqsupseteq j$, we observe that
		\begin{equation} \label{121214_04}
		(j \to i'')' = j.
		\end{equation}
		Hence
		$$
		\begin{array}{lcll}
		i \to j  & = & i \to (j \to i'')' & \mbox{by (\ref{121214_04})} \\
		& = & i \to (j \to i)' & \mbox{} \\
		& = & [(i \to j) \to i]' & \mbox{by Lemma \ref{general_properties2} (\ref{291014_09})} \\
		& = & [(0 \to j) \to i]' & \mbox{by Lemma \ref{general_properties2} (\ref{291014_10})} \\
		& = & [0' \to i]' & \mbox{by Lemma \ref{121214_01} since } j \sqsupseteq 0 \\
		& = & i' & \mbox{by Lemma \ref{general_properties_equiv} (\ref{TXX})} \\
		& = & max(i' , j) & \mbox{since } i' \sqsupseteq j
		\end{array}
		$$

		\item {\bf Case 1.1.2}: $i' < j$.\\
		Since $i' \sqsubseteq j$, we have 
		\begin{equation} \label{121214_03}
		(i' \to j')' = i'.
		\end{equation}
		Therefore,
		$$
		\begin{array}{lcll}
		i \to j & = & i'' \to j & \mbox{} \\
		& = & (i' \to j')'' \to j & \mbox{by (\ref{121214_03})} \\
		& = & (i' \to j') \to j & \mbox{} \\
		& = & (i' \to 0') \to j & \mbox{by Lemma \ref{general_properties2} (\ref{281014_05})} \\
		& = & (0 \to i) \to j & \mbox{by Lemma \ref{general_properties} (\ref{cuasiConmutativeOfImplic2})} \\
		& = & 0' \to j & \mbox{by Lemma \ref{121214_01} since } i \sqsupseteq 0 \\
		& = & j & \mbox{by Lemma \ref{general_properties_equiv} (\ref{TXX})} \\
		& = & max(i' , j) & \mbox{since } i' \sqsubseteq j
		\end{array}
		$$
                   \end{itemize}
		\item {\bf Case 1.2}: 	$i = 0$.\\	
		Using Lemma \ref{121214_01} and Lemma \ref{last_element}, $0 \to j = 0' = max (0' , j)$.
		
		\item {\bf Case 1.3}: $i < 0$.
		
		$$
		\begin{array}{lcll}
	i \to j	& = & (0 \to i) \to j & \mbox{by Lemma \ref{lemma_properties_chains} (\ref{271014_01})} \\
		& = & (i' \to 0') \to j & \mbox{} \\
		& = & (i \to 0') \to j & \mbox{by Lemma \ref{lemma_properties_chains} (\ref{271014_10})} \\
		& = & [(j' \to i) \to (0' \to j)']' & \mbox{by (I)} \\
		& = & [(j' \to i) \to j']' & \mbox{} \\ 
		& = & [(0 \to i) \to j']' & \mbox{by Lemma \ref{general_properties2} (\ref{291014_10})} \\
		& = & (i \to j')' & \mbox{by Lemma \ref{lemma_properties_chains} (\ref{271014_01})} \\
		& = & i & \mbox{since } i \sqsubset j \\
		& = & min(i,  j) & \mbox{} 
		\end{array}
		$$

		\end{itemize}
		
		\item {\bf Case 2}: $j < 0$. \\  
		It is useful to consider the following subcases:
		\begin{itemize}
		
		\item {\bf Case 2.1}: $i > 0$ 
		
		$$
		\begin{array}{lcll}
	i \to j	& = & i \to j' & \mbox{by Lemma \ref{lemma_properties_chains} (\ref{271014_10})} \\
		& = & i \to (j \to i')'' & \mbox{since } j \sqsubset i \\
		& = & i \to (j \to i') & \mbox{} \\
		& = & j \to i' & \mbox{by Lemma \ref{general_properties2} (\ref{281114_01})} \\
		& = & (j \to i')'' & \mbox{} \\
		& = & j' & \mbox{since } j \sqsubset i \\
		& = & j & \mbox{by Lemma \ref{lemma_properties_chains} (\ref{271014_10})} \\
		& = & min(i,  j) & \mbox{} 
		\end{array}
		$$

		\item {\bf Case 2.2}: $i = 0$.


		$$
		\begin{array}{lcll}
	i \to j	& = & 0 \to j & \mbox{} \\
		& = & j & \mbox{by Lemma \ref{lemma_properties_chains} (\ref{271014_01})} \\
		& = & min(i, j) & \mbox{} 
		\end{array}
		$$

		\item {\bf Case 2.3}: $i < 0$.
		       \begin{itemize}
		 \item {\bf Case 2.3.1}: $i \leq j$.\\
		 As $i \sqsubseteq j$, we have 
		\begin{equation} \label{121214_05}
		(i \to j')' = i.
		\end{equation}
		 Observe
		$$
		\begin{array}{lcll}
		i \to j & = & i \to j' & \mbox{by Lemma \ref{lemma_properties_chains} (\ref{271014_10})} \\
		& = & (i \to j')'' & \mbox{} \\
		& = & i' & \mbox{by (\ref{121214_05})} \\
		& = & i & \mbox{by Lemma \ref{lemma_properties_chains} (\ref{271014_10})} \\
		& = & min(i,  j). & \mbox{}
		\end{array}
		$$

		 \item {\bf Case 2.3.2}: $i > j$.
		 We have 
		\begin{equation} \label{121214_06}
		(j \to i')' = j.
		\end{equation}
		as $j \sqsubseteq i$. Hence
		$$
		\begin{array}{lcll}
		i \to j & = & j' \to i' & \mbox{by Lemma \ref{tildeInverImplic}} \\
		& = & j \to i' & \mbox{by Lemma \ref{lemma_properties_chains} (\ref{271014_10})} \\
		& = & (j \to i')'' & \mbox{} \\
		& = & j' & \mbox{by (\ref{121214_06})} \\
		& = & j & \mbox{by Lemma \ref{lemma_properties_chains} (\ref{271014_10})} \\
		& = & min(j , i) & \mbox{}
		\end{array}
		$$
		     \end{itemize}
                 \end{itemize}
		 \item {\bf Case 3}: $j = 0$.
		 \begin{itemize}
		 \item {\bf Case 3.1}: $i \geq 0$.\\
		 By Corollary \ref{121214_08}, as $i \sqsupseteq 0$, we have that $i' = i \to 0 \sqsupseteq 0$. Hence $i \to 0 = i' = max( i',  0)$.

		  \item {\bf Case 3.2}: $i < 0$.
		We have that 
		$$
		\begin{array}{lcll}
	i \to j	& = & i \to 0 & \mbox{} \\
		& = & i' & \mbox{} \\
		& = & i & \mbox{by Lemma \ref{lemma_properties_chains} (\ref{271014_10})} \\
		& = & min(i,  j) & \mbox{} 
		\end{array}
		$$
                   \end{itemize}
					
	\end{itemize}
Hence $\mathbf A \cong \langle [-n;m]; \Rightarrow, 0 \rangle$.
\end{Proof}

\medskip

%
%
The following theorem, our second main result, is now immediate from the preceding results.
\begin{theorem}
	There are $n$ non-isomorphic $I_{2,0}$-chains of size $n$, for $n \in \mathbb{N}$.
\end{theorem}

\appendix
\section{Appendix: Proofs}
We would like to mention here that the identity: $x'' \approx x$ is used in these proofs frequently without explicit mention.\\

\begin{Proof} {\bf of Lemma \ref{general_properties2}}:
Items (\ref{281014_05}) to (\ref{281114_01}) are proved in \cite{CoSa2015aI}. The proofs of (\ref{291014_08}) to (\ref{031214_16}) are given in \cite{CoSa2015semisimple}.
	Let $a,b,c,d \in A$.
	\begin{itemize}
		
		%
		\item[(\ref{281014_06})]
		$$
		\begin{array}{lcll}
		(b \to c) \to a & = & [(b \to c) \to a]' \to [(b \to c) \to a] & \mbox{by Lemma \ref{general_properties_equiv} (\ref{LeftImplicationwithtilde})} \\
		& = & [(a' \to b) \to (c \to a)']'' \to [(b \to c) \to a]  & \mbox{from (I)} \\
		& = & [(a' \to b) \to (c \to a)'] \to [(b \to c) \to a] & \mbox{}
		\end{array}
		$$
		\item[(\ref{201114_02})]
		$$
		\begin{array}{lcll}
		[[0 \to (a \to b)'] \to (0 \to b')']' & = & [\{0 \to (a \to b)'\} \to (b \to 0')']' & \mbox{by Lemma \ref{general_properties} (\ref{cuasiConmutativeOfImplic2})} \\
		& = & [(a \to b)' \to b] \to 0' & \mbox{by (I)} \\
		& = & 0 \to [(a \to b)' \to b]' & \mbox{by Lemma \ref{general_properties} (\ref{cuasiConmutativeOfImplic2})} \\
		& = & (a \to b) \to (0 \to b') & \mbox{by (\ref{071114_01})} \\
		& = & 0 \to [(a \to b) \to b'] & \mbox{by (\ref{071114_04})} \\
		& = & 0 \to [(a \to 0') \to b'] & \mbox{by  (\ref{281014_05})} \\
		& = & (a \to 0') \to (0 \to b') & \mbox{by (\ref{071114_04})} \\
		& = & (0 \to a') \to (0 \to b') & \mbox{by Lemma \ref{general_properties} (\ref{cuasiConmutativeOfImplic2})} \\
		& = & a' \to (0 \to b') & \mbox{by (\ref{311014_06})} \\
		& = & [(0 \to a) \to (0 \to b')']' & \mbox{by (\ref{031114_06})} \\
		& = & [(0 \to a) \to (b \to 0')']' & \mbox{by Lemma \ref{general_properties} (\ref{cuasiConmutativeOfImplic2})} \\
		& = & (a \to b) \to 0' & \mbox{by (I)} \\
		& = & 0 \to (a \to b)' & \mbox{by Lemma \ref{general_properties} (\ref{cuasiConmutativeOfImplic2})}
		\end{array}
		$$
		\item[(\ref{201114_03})]
		$$
		\begin{array}{lcll}
		0 \to [(a \to b) \to c]   & = & 0 \to [(c' \to a) \to (b \to c)']' & \mbox{by (I)} \\
		& \sqsubseteq & 0 \to (b \to c)'' & \mbox{by (\ref{201114_02})} \\
		& = & 0 \to (b \to c) & \mbox{}
		\end{array}
		$$
		\item[(\ref{201114_08})]
		$$
		\begin{array}{lcll}
		a' \to (b \to 0')' & = & (a \to 0) \to (b \to 0')' & \mbox{} \\
		& = & [\{(b \to 0') \to a\} \to \{0 \to (b \to 0')'\}']' & \mbox{by (I)} \\
		& = & [\{(b \to 0') \to a\} \to \{0 \to (0 \to b')'\}']' & \mbox{by Lemma \ref{general_properties} (\ref{cuasiConmutativeOfImplic2})} \\
		& = & [\{(b \to 0') \to a\} \to (0 \to b)']' & \mbox{by (\ref{031114_07})} \\
		& = & [\{(0 \to b') \to a\} \to (0 \to b)']' & \mbox{by Lemma \ref{general_properties} (\ref{cuasiConmutativeOfImplic2})} \\
		& = & [[\{0 \to (0 \to b)'\} \to a] \to (0 \to b)']' & \mbox{by (\ref{031114_07})} \\
		& = & [a \to (0 \to b)']' & \mbox{by (\ref{291014_08})}
		\end{array}
		$$
		\item[(\ref{251114_03})]
		$$
		\begin{array}{lcll}
		[(0 \to a) \to b]' & = & [(b \to a) \to b]' & \mbox{by (\ref{291014_10})} \\
		& = & [b \to (a \to b)']'' & \mbox{by (\ref{291014_09})} \\
		& = & b \to (a \to b)' & \mbox{}
		\end{array}
		$$
		\item[(\ref{271114_01})]
		
		$$
		\begin{array}{lcll}
		[a \to (b \to 0')']' & = & [a \to (0 \to b')']' & \mbox{by Lemma \ref{general_properties} (\ref{cuasiConmutativeOfImplic2})} \\
		& = & a' \to (b' \to 0')' & \mbox{by (\ref{201114_08})} \\
		& = & a' \to (0 \to b)' & \mbox{by Lemma \ref{general_properties} (\ref{cuasiConmutativeOfImplic2})}
		\end{array}
		$$
		\item[(\ref{271114_02})]
		$$
		\begin{array}{lcll}
		b' \to a' & = & (b \to 0) \to a' & \mbox{} \\
		& = & [(a \to b) \to (0 \to a')']' & \mbox{by (I)} \\
		& = & [(a \to b) \to (a \to 0')']' & \mbox{by Lemma \ref{general_properties} (\ref{cuasiConmutativeOfImplic2})} \\
		& = & (a \to b)' \to (0 \to a)' & \mbox{by (\ref{271114_01}) with } x = a \to b, y = a
		\end{array}
		$$
		\item[(\ref{011214_01})]
		
		$$
		\begin{array}{lcll}
		(0 \to a)' \to (0 \to b)' & = & [(0 \to a) \to 0] \to (0 \to b)' \\
		& = & [\{(0 \to b) \to (0 \to a)\} \to \{0 \to (0 \to b)'\}']' &\\
		&    & \hspace{4cm}       \mbox{by (I)} \\
		& = & [\{0 \to (0 \to b)'\} \to (0 \to b)] \to [(0 \to a) \to \{0 \to (0 \to b)'\}']' & \\
		&  &  \hspace{4cm} \mbox{by (I)} \\
		& = & [\{(0 \to b) \to (0 \to b)'\} \to (0 \to b)] \to [(0 \to a) \to \{0 \to (0 \to b)'\}']' &\\
		&   &   \hspace{4cm} \mbox{by (\ref{291014_10})} \\
		& = & [(0 \to b)' \to (0 \to b)] \to [(0 \to a) \to \{0 \to (0 \to b)'\}']' &  \\
		&  &   \hspace{4cm}      \mbox{by Lemma \ref{general_properties_equiv} (\ref{LeftImplicationwithtilde})} &  \\
		& = &  (0 \to b) \to [(0 \to a) \to \{0 \to (0 \to b)'\}']' &  \\
		&  &  \hspace{4cm}    \mbox{by Lemma \ref{general_properties_equiv} (\ref{LeftImplicationwithtilde})} &  \\
		& = & (0 \to b) \to [(0 \to a) \to (0 \to b')']' & \\
		&  &  \hspace{4cm} \mbox{by (\ref{031114_07})} \\
		& = & (0 \to b) \to [(0 \to a) \to (b \to 0')']' &  \\
		&  &  \hspace{4cm} \mbox{by Lemma \ref{general_properties} (\ref{cuasiConmutativeOfImplic2})} &  \\
		& = & (0 \to b) \to [(a \to b) \to 0'] & \\ 
		&  &  \hspace{4cm} \mbox{by (I)} \\
		& = & (0 \to b) \to [0 \to (a \to b)'] &  \\
		&  &   \hspace{4cm} \mbox{by Lemma \ref{general_properties} (\ref{cuasiConmutativeOfImplic2})} &  \\
		& = & 0 \to [(0 \to b) \to (a \to b)'] & \\
		&  &  \hspace{4cm} \mbox{by (\ref{071114_04})} \\
		& = & 0 \to (a \to b)' & \\ 
		&  &  \hspace{4cm} \mbox{by (\ref{291014_06})} \\
		& = & (a \to b) \to 0'  &  \\
		&  &  \hspace{4cm} \mbox{by Lemma \ref{general_properties} (\ref{cuasiConmutativeOfImplic2})} &  \\
		& = & [(0 \to a) \to (b \to 0')']' & \\ 
		&  &  \hspace{4cm}\mbox{by (I)} \\
		& = & [(0 \to a) \to (0 \to b')']' &  \\
		&  &  \hspace{4cm} \mbox{by Lemma \ref{general_properties} (\ref{cuasiConmutativeOfImplic2})} &  \\
		& = & [a' \to (0 \to b')]'' & \\
		&  &  \hspace{4cm} \mbox{by (\ref{031114_06})} \\
		& = & a' \to (0 \to b') & \mbox{} \\
		& = & 0 \to (a' \to b') &  \hspace{-7cm} \mbox{by (\ref{071114_04})}.
		\end{array}
		$$
		\item[(\ref{011214_03})]
		$$
		\begin{array}{lcll}
		[(a \to b)' \to \{b \to (a \to b)'\}']' & = & [(a \to b)' \to b] \to (a \to b)' & \mbox{by (\ref{291014_09})} \\
		& = & (0 \to b) \to (a \to b)' & \mbox{by (\ref{291014_10})} \\
		& = & (a \to b)' & \mbox{by (\ref{291014_06})}
		\end{array}
		$$
		\item[(\ref{011214_04})]
		
		%
		$$
		\begin{array}{lcll}
		(0 \to a) \to b & = & (b \to a) \to b & \mbox{by (\ref{291014_10})} \\
		& = & [b \to (a \to b)']' & \mbox{by (\ref{291014_09})} \\
		& \sqsubseteq & (a \to b)' \to [b \to (a \to b)']' & \mbox{by (\ref{011214_03}) with } x = b, y = (a \to b)' \\
		& = & [\{(a \to b)' \to b\} \to (a \to b)']' & \mbox{by (\ref{291014_09})} \\
		& = & [(0 \to b) \to (a \to b)']' & \mbox{by (\ref{291014_10})} \\
		& = & (a \to b)'' & \mbox{by (\ref{291014_06})} \\
		& = & a \to b & \mbox{ since } x'' \approx x 
		\end{array}
		$$
		\item[(\ref{031214_07})]
		$$
		\begin{array}{lcll}
		[\{a \to (b \to a)'\} \to a'']'  & = & [\{a \to (b \to a)'\} \to a]' & \mbox{} \\
		& = & [\{0 \to (b \to a)'\} \to a]' & \mbox{by (\ref{291014_10})} \\
		& = & [\{(b \to a) \to 0'\} \to a]' & \mbox{by Lemma \ref{general_properties} (\ref{cuasiConmutativeOfImplic2})} \\
		& = & [\{(b \to a) \to a'\} \to a]' & \mbox{by (\ref{281014_05})} \\
		& = & [\{(b \to 0') \to a'\} \to a]' & \mbox{by (\ref{281014_05})} \\
		& = & [\{(b \to 0') \to 0'\} \to a]' & \mbox{by (\ref{281014_05})} \\
		& = & [\{(b \to 0'') \to 0'\} \to a]' & \mbox{by (\ref{281014_05})} \\
		& = & [\{(b \to 0) \to 0'\} \to a]' & \mbox{} \\
		& = & [(b' \to 0') \to a]' & \mbox{} \\
		& = & [(0 \to b) \to a]' & \mbox{by Lemma \ref{general_properties} (\ref{cuasiConmutativeOfImplic2})} \\
		& = & [(a \to b) \to a]' & \mbox{by (\ref{291014_10})} \\
		& = & a \to (b \to a)' & \mbox{by (\ref{291014_09})}
		\end{array}
		$$
		%
	\end{itemize}
\end{Proof}

\begin{Proof} {\bf of Lemma \ref{pretransitivity1}}
	\begin{itemize}
		\item[(\ref{031114_10})]
		Observe that by Lemma \ref{general_properties} (\ref{cuasiConmutativeOfImplic2}), Lemma \ref{general_properties2} (\ref{281014_05}) and the hypothesis we have that $(0 \to a') \to b = (a \to 0') \to b = (a \to b') \to b = (a \to b')'' \to b = a' \to b$. 
		
		\item[(\ref{031114_11})]
		$$
		\begin{array}{lcll}
		b \to a' & = & [(0 \to a') \to b] \to a' & \mbox{by Lemma \ref{general_properties2} (\ref{291014_08})} \\
		& = & (a' \to b) \to a' & \mbox{from (\ref{031114_10})} \\
		& = & (0 \to b) \to a' & \mbox{by Lemma \ref{general_properties2} (\ref{291014_10})}.
		\end{array}
		$$

		\item[(\ref{031114_12})]
		
		$$
		\begin{array}{lcll}
		b \to a' & = & (0 \to b) \to a' & \mbox{from (\ref{031114_11})} \\
		& = & (0 \to b) \to (a \to b')'' & \mbox{by hypothesis} \\
		& = & (0 \to b) \to (a \to b') & \mbox{} \\
		& = & (0 \to b'') \to (a \to b') & \mbox{} \\
		& = & a \to b' & \mbox{by Lemma \ref{general_properties2} (\ref{281014_07})} \\
		& = & (a \to b')'' & \mbox{} \\
		& = & a' & \mbox{by hypothesis}
		\end{array}
		$$
		
		\item[(\ref{201114_04})]
		$$
		\begin{array}{lcll}
		0 \to (a' \to b) & = & a' \to (0 \to b) & \mbox{by Lemma \ref{general_properties2} (\ref{071114_04})} \\
		& = & 0 \to (a \to b')' & \mbox{by Lemma \ref{general_properties2} (\ref{071114_01})} \\
		& = & 0 \to a & \mbox{by hypothesis}
		\end{array}
		$$
		
		\item[(\ref{211114_03})]
		By hypothesis and (I) we have that $(d \to a) \to b' = [(b \to d) \to (a \to b')']' = [(b \to d) \to a]'$.

		\item[(\ref{211114_01})]
		$$
		\begin{array}{lcll}
		[\{d \to (0 \to b')\} \to a]' & = & (a' \to d) \to [(0 \to b') \to a]' & \mbox{by (I)} \\
		& = & (a' \to d) \to [(a \to b') \to a]' & \mbox{by Lemma \ref{general_properties2} (\ref{291014_10})} \\
		& = & (a' \to d) \to [(a \to b')'' \to a]' & \mbox{} \\
		& = & (a' \to d) \to (a' \to a)' & \mbox{by hypothesis} \\
		& = & (a' \to d) \to a' & \mbox{by Lemma \ref{general_properties_equiv} (\ref{LeftImplicationwithtilde})} \\
		& = & (a' \to d) \to (0' \to a)' & \mbox{by Lemma \ref{general_properties_equiv} (\ref{TXX})} \\
		& = & [(d \to 0') \to a]' & \mbox{by (I)} \\
		& = & (0 \to d) \to a' & \mbox{by Lemma \ref{general_properties2} (\ref{071114_02})}
		\end{array}
		$$
		
		\item[(\ref{211114_02})]
		
		$$
		\begin{array}{lcll}
		a \to [(a' \to d) \to \{(0 \to a) \to b'\}] & = & a \to [(a' \to d) \to \{(b \to 0) \to (a \to b')'\}'] & \mbox{} \\
		&  & \mbox{by (I)} & \mbox{} \\
		& = & a \to [(a' \to d) \to \{(b \to 0) \to a\}'] & \mbox{} \\
		&  & \mbox{by hypothesis} & \mbox{} \\
		& = & a \to [\{d \to (b \to 0)\} \to a]' & \mbox{} \\
		&  & \mbox{by (I)} & \mbox{} \\
		& = & \left[[d \to (0 \to (b \to 0))] \to a\right]' & \mbox{} \\
		&  & \mbox{by Lemma \ref{general_properties2} (\ref{171114_01}) with } x = d, y = b \to 0, z = a & \mbox{} \\
		& = & \left[[d \to (0 \to b')] \to a\right]' & \mbox{} \\
		& = & (0 \to d) \to a' & \mbox{} \\
		&  & \mbox{by (\ref{211114_01})} & \mbox{} 
		\end{array}
		$$
		
		\item[(\ref{211114_04})]
		$$
		\begin{array}{lcll}
		a \to ((d \to a) \to b') & = & a \to [(b \to d) \to a]' & \mbox{by (\ref{211114_03})} \\
		& = & a'' \to [(b \to d) \to a]' & \mbox{} \\
		& = & (a' \to 0) \to [(b \to d) \to a]' & \mbox{} \\
		& = & [\{0 \to (b \to d)\} \to a]' & \mbox{by (I)} \\
		& = & [\{(b \to d)' \to 0'\} \to a]' & \mbox{} \\
		& = & [\{((b \to d) \to 0) \to 0'\} \to a]' & \mbox{} \\
		& = & [\{((b \to d) \to 0') \to 0'\} \to a]' & \mbox{by Lemma \ref{general_properties2} (\ref{281014_05})} \\
		& = & [\{((b \to d) \to 0') \to a'\} \to a]' & \mbox{by Lemma \ref{general_properties2} (\ref{281014_05})} \\
		& = & [\{((b \to d) \to a) \to a'\} \to a]' & \mbox{by Lemma \ref{general_properties2} (\ref{281014_05})} \\
		& = & [\{((b \to d) \to a) \to 0'\} \to a]' & \mbox{by Lemma \ref{general_properties2} (\ref{281014_05})} \\
		& = & [\{0 \to ((b \to d) \to a)'\} \to a]' & \mbox{} \\
		& = & [\{0 \to ((d \to a) \to b')\} \to a]' & \mbox{by (\ref{211114_03})} \\
		& = & [\{a \to ((d \to a) \to b')\} \to a]' & \mbox{by Lemma \ref{general_properties2} (\ref{291014_10})} \\
		& = & (a' \to a) \to [\{(d \to a) \to b'\} \to a]' & \mbox{by (I)} \\
		& = & a \to [\{(d \to a) \to b'\} \to a]' & \mbox{by Lemma \ref{general_properties_equiv} (\ref{LeftImplicationwithtilde})} \\
		& = & a \to [\{a' \to (d \to a)\} \to (b' \to a)'] & \mbox{by (I)} \\
		& = & a \to [\{a' \to (d \to a)\} \to \{(b \to 0) \to a\}'] & \mbox{} \\
		& = & a \to [\{a' \to (d \to a)\} \to \{(b \to 0) \to (a \to b')'\}'] & \mbox{by hypothesis} \\
		& = & a \to [\{a' \to (d \to a)\} \to \{(0 \to a) \to b'\}] & \mbox{by (I)}  \\
		& = & [0 \to (d \to a)] \to a' & \mbox{by (\ref{211114_02}) with } d := d \to a \\
		& = & a \to (d \to a)' & \mbox{by Lemma \ref{general_properties2} (\ref{181114_11})}
		\end{array}
		$$
		
		\item[(\ref{211114_06})]
		$$
		\begin{array}{lcll}
		[0 \to (b \to d)] \to a  & = & [(a' \to 0) \to ((b \to d) \to a)']' & \mbox{by (I)} \\
		& = & [a \to ((b \to d) \to a)']' & \mbox{} \\
		& = & [a \to ((d \to a) \to b')]' & \mbox{by (\ref{211114_03})} \\
		& = & [a \to (d \to a)']' & \mbox{by (\ref{211114_04})} \\
		& = & (a \to d) \to a & \mbox{by Lemma \ref{general_properties2} (\ref{291014_09})} \\
		& = & (0 \to d) \to a & \mbox{by Lemma \ref{general_properties2} (\ref{291014_10})}
		\end{array}
		$$
		
		\item[(\ref{211114_07})]
		$$
		\begin{array}{lcll}
		(b \to (a \to d)) \to a & = & [(a' \to b) \to \{(a \to d) \to a\}']' & \mbox{by (I)} \\
		& = & [(a' \to b) \to \{(0 \to d) \to a\}']' & \mbox{by Lemma \ref{general_properties2} (\ref{291014_10})} \\
		& = & [b \to (0 \to d)] \to a & \mbox{by (I)} \\
		& = & [0 \to (b \to d)] \to a & \mbox{by Lemma \ref{general_properties2} (\ref{071114_04})} \\
		& = & (0 \to d) \to a & \mbox{by (\ref{211114_06})}
		\end{array}
		$$

		\item[(\ref{031214_11})]
		$$
		\begin{array}{lcll}
		b \to (0 \to a')  & = & (0 \to b) \to (0 \to a') & \mbox{by Lemma \ref{general_properties2} (\ref{311014_06})} \\
		& = & 0 \to [(0 \to b) \to a'] & \mbox{by Lemma \ref{general_properties2} (\ref{071114_04})}  \\
		& = & 0 \to [(a \to 0) \to (b \to a')']' & \mbox{by (I)} \\
		& = & 0 \to [a' \to (b \to a')']' & \mbox{} \\
		& = & 0 \to (a' \to a'')' & \mbox{by (\ref{031114_12})} \\
		& = & 0 \to (a' \to a)' & \mbox{} \\
		& = & 0 \to a' & \mbox{by Lemma \ref{general_properties_equiv} (\ref{LeftImplicationwithtilde})}
		\end{array}
		$$
		
		\item[(\ref{031214_13})]
		From (I) and by hypothesis we have that $[(d \to a) \to b']' = (b \to d) \to (a \to b')' = (b \to d) \to a$.
		
		\item[(\ref{031214_14})]
		$$
		\begin{array}{lcll}
		a' \to b & = & (a \to b') \to b & \mbox{by hypothesis} \\
		& = & [(b' \to a) \to (b' \to b)']' & \mbox{using (I)} \\
		& = & [(b' \to a) \to b']' & \mbox{by Lemma \ref{general_properties_equiv} (\ref{LeftImplicationwithtilde})} \\
		& = & b' \to (a \to b')' & \mbox{by Lemma \ref{general_properties2} (\ref{291014_09})} \\
		& = & b' \to a & \mbox{by hypothesis}
		\end{array}
		$$
		
		\item[(\ref{031214_15})]
		$$
		\begin{array}{lcll}
		(d \to 0') \to (a' \to b) & = & (a' \to b)' \to  [(d \to 0') \to (a' \to b)] & \\
		&  &  \hspace{4cm} \mbox{by Lemma \ref{general_properties2} (\ref{281114_01})} \\
		& = & (a' \to b)' \to  [(d \to 0') \to (b' \to a)] & \\
		&  &   \hspace{4cm}\mbox{by (\ref{031214_14})} \\
		& = & (a' \to b)' \to  [(d \to 0') \to \{(0 \to a) \to b'\}'] & \\ 
		&  & \hspace{4cm}\mbox{by (\ref{031214_13}) with } d = 0 \\
		& = & (a' \to b)' \to  [\{(0 \to a) \to b'\} \to \{d \to ((0 \to a) \to b')\}'] &\\
		&  &  \hspace{4cm} \mbox{by Lemma \ref{general_properties2} (\ref{031214_16}) with }  \\
		&  &   \hspace{4cm} x = (0 \to a) \to b', y = d \\
		& = & (a' \to b)' \to  [(b' \to a)' \to \{d \to ((0 \to a) \to b')\}'] &\\
		&  &  \hspace{4cm} \mbox{by (\ref{031214_13}) with } d = 0 \\
		& = & (a' \to b)' \to  [(a' \to b)' \to \{d \to ((0 \to a) \to b')\}'] & \\
		&  &  \hspace{4cm} \mbox{by (\ref{031214_14})} \\
		& = & (a' \to b)' \to [d \to \{(0 \to a) \to b'\}]' & \\
		&  &  \hspace{4cm} \mbox{by Lemma \ref{general_properties2} (\ref{031114_04})} \\
		& = & (b' \to a)' \to [d \to \{(0 \to a) \to b'\}]' & \\
		&  &  \hspace{4cm} \mbox{by (\ref{031214_14})} \\
		& = & [(0 \to a) \to b'] \to [d \to \{(0 \to a) \to b'\}]' & \\
		&  &  \hspace{4cm}\mbox{by (\ref{031214_13}) with } d = 0 \\
		& = & [(0 \to d) \to \{(0 \to a) \to b'\}]' &\\
		&  &  \hspace{4cm} \mbox{by Lemma \ref{general_properties2} (\ref{251114_03}) with }  \\
		&  &   \hspace{4cm} x = d, y = (0 \to a) \to b' \\
		& = & [(b' \to d) \to \{(0 \to a) \to b'\}]' & \\
		&  &  \hspace{4cm} \mbox{by Lemma \ref{general_properties2} (\ref{181114_07}) with } \\
		&  &   \hspace{4cm} x = b, y = d, z = a \\
		& = & [(b' \to d) \to (b' \to a)']'  & \\
		&  &  \hspace{4cm}\mbox{by (\ref{031214_13}) with } d = 0 \\
		& = & [(b' \to d) \to (a' \to b)']' &\\
		&  &  \hspace{4cm} \mbox{by (\ref{031214_14})} \\
		& = & (d \to a') \to b &  \hspace{-5.1cm}\mbox{by (I)}.
		\end{array}
		$$
		
		\item[(\ref{031214_21})]
		$$
		\begin{array}{lcll}
		[(0 \to a') \to b]' & = & [(a \to 0') \to b]' & \mbox{by Lemma \ref{general_properties} (\ref{cuasiConmutativeOfImplic2})} \\
		& = & (0 \to a) \to b' & \mbox{by Lemma \ref{general_properties2} (\ref{071114_02})}
		\end{array}
		$$
		
		\item[(\ref{031214_17})]
		$$
		\begin{array}{lcll}
		(a' \to b)' & = & [(0 \to a') \to b]' & \mbox{by (\ref{031114_10})} \\
		& = & (0 \to a) \to b' & \mbox{by (\ref{031214_21})}
		\end{array}
		$$

		\item[(\ref{031214_19})]
		$$
		\begin{array}{lcll}
		[\{b' \to ((b \to d) \to a)\} \to (0 \to b)']'   & = & [\{b' \to ((d \to a) \to b')'\} \to (0 \to b)']' &\\
		&  &  \hspace{4cm} \mbox{by (\ref{031214_13})} \\
		& = & [\{(d \to a) \to b'\}' \to 0] \to b & \\
		&  &  \hspace{4cm} \mbox{by (I)} \\
		& = & [(d \to a) \to b'] \to b & \mbox{} \\
		& = & b' \to [(d \to a) \to b']' & \\
		&  &  \hspace{4cm} \mbox{by Lemma \ref{general_properties2} (\ref{071114_05}) with } \\
		&  &   \hspace{4cm}   x = d \to a, y = b' \\
		& = & b' \to [(b \to d) \to a] &   \hspace{-4.5cm}\mbox{by (\ref{031214_13})}.
		\end{array}
		$$

	\end{itemize}
	
\end{Proof}

\end{document}